# LARGE DEVIATIONS FOR RANDOM WALKS UNDER SUBEXPONENTIALITY: THE BIG-JUMP DOMAIN

By D. Denisov,[1] A. B. Dieker [2] and V. Shneer[1]

*Heriot-Watt University, IBM Research and EURANDOM*

For a given one-dimensional random walk $\{S_n\}$ with a subexponential step-size distribution, we present a unifying theory to study the sequences $\{x_n\}$ for which $\mathsf{P}\{S_n > x\} \sim n\mathsf{P}\{S_1 > x\}$ as $n \to \infty$ uniformly for $x \geq x_n$. We also investigate the stronger "local" analogue, $\mathsf{P}\{S_n \in (x, x+T]\} \sim n\mathsf{P}\{S_1 \in (x, x+T]\}$. Our theory is self-contained and fits well within classical results on domains of (partial) attraction and local limit theory.

When specialized to the most important subclasses of subexponential distributions that have been studied in the literature, we reproduce known theorems and we supplement them with new results.

**1. Introduction.** In general, it poses a challenge to find the exact asymptotics for probabilities that tend to zero. However, due to the vast set of available tools, a great deal is known about probabilities arising from a one-dimensional random walk $\{S_n\}$. For instance, under Cramér's condition on the step-size distribution, the famous Bahadur–Ranga Rao theorem describes the deviations of $S_n/n$ from its mean; see, for instance, Höglund [22]. Other random walks with well-studied (large) deviation behavior include those with step-size distributions for which Cramér's condition *does not* hold.

*Large deviations under subexponentiality.* The present paper studies large deviations for random walks with *subexponential* step-size distributions on the real line. These constitute a large class of remarkably tractable distributions for which Cramér's condition does not hold. The resulting random

Received March 2007; revised November 2007.
[1]Supported by the Dutch BSIK project (BRICKS) and the EURO-NGI project.
[2]Supported by the Science Foundation Ireland Grant SFI04/RP1/I512, and by the Netherlands Organization for Scientific Research (NWO) under Grant 631.000.002.

walks have the property that there exists some sequence $\{x_n\}$ (depending on the step-size distribution) for which [9]

$$\lim_{n\to\infty} \sup_{x\geq x_n} \left| \frac{\mathsf{P}\{S_n > x\}}{n\mathsf{P}\{S_1 > x\}} - 1 \right| = 0. \tag{1}$$

The intuition behind the factor $n$ is that a single big increment causes $S_n$ to become large, and that this "jump" may occur at each of the $n$ epochs. Given a subexponential step-size distribution, it is our aim to characterize sequences $\{x_n\}$ for which (1) holds. In other words, we are interested in (the boundary of) the big-jump domain.

The big-jump domain has been well studied for special classes of subexponential distributions. Overviews are given in Embrechts, Klüppelberg and Mikosch [14], Section 8.6, S. Nagaev [33] and Mikosch and A. Nagaev [30]. Due to its importance in applications (e.g., [10]), there is a continuing interest in this topic. Work published after 2003 includes Baltrūnas, Daley and Klüppelberg [2], Borovkov and Mogulśkiĭ [7], Hult et al. [23], Jelenković and Momčilović [25], Konstantinides and Mikosch [27], Ng et al. [35] and Tang [44]. Finally, we also mention the important articles by Pinelis [38, 39] and Rozovskii [40, 41]. Pinelis studies large deviations for random walks in Banach spaces, while Rozovskii investigates general deviations from the mean, *beyond* the big-jump domain. Our paper owes much to Rozovskii's work.

*Novelties.* Although the sequences for which (1) holds have been characterized for certain subclasses of subexponential distributions, the novelty of our work is twofold:

- we present a *unified* theory within the framework of subexponentiality, which fits well within classical results on domains of (partial) attraction and local limit theory, and
- we also study the local analogue of (1); that is, for a given $T > 0$, we study the $x$-domain for which $\mathsf{P}\{S_n \in (x, x+T]\}$ is uniformly approximated by $n\mathsf{P}\{S_1 \in (x, x+T]\}$.

When specialized to the classes of subexponential distributions studied in the literature, our theory reproduces the sharpest known results with short proofs. Moreover, in some cases it allows to improve upon the best-known boundaries by several orders of magnitude, as well as to derive entirely new results.

By presenting a unified large-deviation theory for subexponential distributions in the big-jump domain, we reveal two effects which play an equally important role. The first effect ensures that having *many "small" steps is unlikely* to lead to the rare event $\{S_n > x\}$, and the second effect requires that the step-size distribution be *insensitive* to shifts on the scale of fluctuations of $S_n$; the latter is known to play a role in the finite-variance case [25, 31]. Since one of these effects typically dominates, this explains the inherently different nature of some of the big-jump boundaries found in the literature.

It is instructive to see how these two effects heuristically solve the large-deviation problem for centered subexponential distributions with unit variance. In this context, the many-small-steps-effect requires that $x \geq J_n$, where $J_n$ satisfies $J_n^2 \sim -2n \log[n\mathsf{P}\{S_1 > J_n\}]$ as $n \to \infty$ [here $f(x) \sim g(x)$ stands for $\lim_x f(x)/g(x) = 1$]. In fact, $J_n$ usually needs to be chosen slightly larger. On the other hand, the insensitivity-effect requires that $x \geq I_n$, where $I_n$ satisfies $\mathsf{P}\{S_1 > I_n - \sqrt{n}\} \sim \mathsf{P}\{S_1 > I_n\}$. After overcoming some technicalities, our theory allows us to show that (1) holds for $x_n = I_n + J_n$. We stress, however, that not only do our results apply to the finite-variance case, but that seemingly "exotic" step-size distributions with infinite mean fit seamlessly into the framework.

The second novelty of our work, the investigation of local asymptotics, also has far-reaching consequences. A significant amount of additional arguments are needed to prove our results in the local case, but local large-deviation theorems are much stronger than their global counterparts. Let us illustrate this by showing that our local results under subexponentiality immediately yield interesting and new theorems within the context of *light* tails. Indeed, given $\gamma > 0$ and a subexponential distribution function $F$ for which $L(\gamma) = \int e^{-\gamma y} F(dy) < \infty$, consider the random walk under the measure $\mathsf{P}^*$ determined by

$$\mathsf{P}^*\{S_1 \in dx\} = \frac{e^{-\gamma x} F(dx)}{\int_{\mathbb{R}} e^{-\gamma y} F(dy)}.$$

Distributions of this form belong to the class which is usually called $\mathcal{S}(\gamma)$ (but $\mathcal{S}(\gamma)$ is larger; see [13]). Suppose that for any $T > 0$, we have $\mathsf{P}\{S_n \in (x, x+T]\} \sim n\mathsf{P}\{S_1 \in (x, x+T]\}$ uniformly for $x \geq x_n$, where $\{S_n\}$ is a $\mathsf{P}$-random walk with step-size distribution $F$ and $\{x_n\}$ does not depend on $T$. Using our local large-deviation results and an elementary approximation argument, we readily obtain that

$$\lim_{n\to\infty} \sup_{x \geq x_n} \left| \frac{\mathsf{P}^*\{S_n > x\}}{nL(\gamma)^{1-n}\mathsf{P}^*\{S_1 > x\}} - 1 \right| = 0.$$

Apart from the one-dimensional random-walk setting, our techniques seem to be suitable to deal with a variety of problems outside the scope of the present paper. For instance, our arguments may unify the results on large deviations for multidimensional random walks [4, 23, 32]. Stochastic recurrences form another challenging area; see [27].

*Outline.* This paper is organized as follows. In Section 2, we introduce four sequences that facilitate our analysis. We also state our main result and outline the idea of the proof. Sections 3–5 contain the proofs of the claims made in Section 2. Two sequences are typically hardest to find, and we derive a series of useful tools to find these sequences in Sections 6 and 7. As a corollary, we obtain a large-deviation result which allows one to conclude that (1) holds with $x_n = an$ for some $a > 0$. In Sections 8 and 9, we work out the most important special cases of our theory. An Appendix treats some notions used in the body of the paper. Appendix A focuses on Karamata theory, while Appendix B discusses the class of subexponential densities.

**2. Main result and the idea of the proof.** We first introduce some notation. Throughout, we study the random walk $\{S_n \equiv \xi_1 + \cdots + \xi_n\}$ with generic step $\xi$. Let $F$ be the step-size distribution, that is, the distribution of $\xi$. We also fix some $T \in (0,\infty]$, and write $F(x+\Delta)$ for $\mathsf{P}\{x < \xi \le x + T\}$, which is interpreted as $\overline{F}(x) \equiv \mathsf{P}\{\xi > x\}$ if $T = \infty$. Apart from these notions, a crucial role in the present paper is also played by $\overline{G}(x) \equiv \mathsf{P}\{|\xi| > x\}$, and the truncated moments $\mu_1(x) \equiv \int_{|y|\le x} yF(dy)$ and $\mu_2(x) \equiv \int_{|y|\le x} y^2 F(dy)$.

We say that $F$ is (locally) long-tailed, written as $F \in \mathcal{L}_\Delta$, if $F(x+\Delta) > 0$ for sufficiently large $x$ and $F(x+y+\Delta) \sim F(x+\Delta)$ for all $y \in \mathbb{R}$. Since this implies that $x \mapsto F(\log x + \Delta)$ is slowly varying, the convergence holds locally uniformly in $y$. The distribution $F$ is (locally) subexponential, written as $F \in \mathcal{S}_\Delta$, if $F \in \mathcal{L}_\Delta$ and $F^{(2)}(x+\Delta) \sim 2F(x+\Delta)$ as $x \to \infty$. Here $F^{(2)}$ is the twofold convolution of $F$. In the local case, for $F$ supported on $[0,\infty)$, the class $\mathcal{S}_\Delta$ has been introduced by Asmussen, Foss and Korshunov [1].

Throughout, both $f(x) \ll g(x)$ and $f(x) = o(g(x))$ as $X \to \infty$ are shorthand for $\lim_{x\to\infty} f(x)/g(x) = 0$, while $f(x) \gg g(x)$ stands for $g(x) \ll f(x)$. We write $f(x) = O(g(x))$ if $\limsup_{x\to\infty} f(x)/g(x) < \infty$, and $f(x) \asymp g(x)$ if $f(x) = O(g(x))$ and $g(x) = O(f(x))$.

With the only exception of our main theorem, Theorem 2.1, all proofs for this section are deferred to Section 3. The proof of Theorem 2.1 is given in Section 4 (global case) and Section 5 (local case).

2.1. *Four sequences; main result.* Our approach relies on four sequences associated to $F$.

*Natural scale.* We say that a sequence $\{b_n\}$ is a *natural-scale sequence* if $\{S_n/b_n\}$ is tight. Recall that this means that for any $\epsilon > 0$, there is some $K > 0$ such that $\mathsf{P}\{S_n/b_n \in [-K,K]\} > 1 - \epsilon$ for all $n$. An equivalent definition is that any subsequence contains a subsequence which converges in distribution. Hence, if $S_n/b_n$ converges in distribution, then $\{b_n\}$ is a natural-scale sequence. For instance, if $\mathsf{E}\{\xi\} = 0$ and $\mathsf{E}\{\xi^2\} < \infty$, then $b \equiv \{\sqrt{n}\}$ is a natural-scale sequence by the central limit theorem.

Due to their prominent role in relation to domain of partial attractions, natural-scale sequences have been widely studied and are well understood; necessary and sufficient conditions for $\{b_n\}$ to be a natural-scale sequence can be found in Section IX.7 of Feller [17]. We stress, however, that we allow for the possibility that $S_n/b_n$ converges in distribution to a degenerate limit; this is typically ruled out in much of the literature. To give an example, suppose that $\mathsf{E}\{\xi\} = 0$ and that $\mathsf{E}\{|\xi|^r\} < \infty$ for some $r \in [1, 2)$. Then $b \equiv \{n^{1/r}\}$ is a natural-scale sequence since $S_n/n^{1/r}$ converges to zero by the Kolmogorov–Marcinkiewicz–Zygmund law of large numbers.

We now collect some facts on natural-scale sequences. First, by the lemma in Section IX.7 of [17] (see also Jain and Orey [24]), we have

$$\lim_{K\to\infty} \sup_n n\overline{G}(Kb_n) = 0 \tag{2}$$

for any natural-scale sequence. The next exponential bound lies at the heart of the present paper.

Lemma 2.1. *For any natural-scale sequence $\{b_n\}$, there exists a constant $C \in (0, \infty)$ such that for any $n \geq 1$, $c \geq 1$ and $x \geq 0$,*

$$\mathsf{P}\{S_n > x, \xi_1 \leq cb_n, \ldots, \xi_n \leq cb_n\} \leq C \exp\left\{-\frac{x}{cb_n}\right\}$$

*and*

$$\mathsf{P}\{|S_n| > x, |\xi_1| \leq cb_n, \ldots, |\xi_n| \leq cb_n\} \leq C \exp\left\{-\frac{x}{cb_n}\right\}.$$

*Insensitivity.* Given a sequence $b \equiv \{b_n\}$, we say that $\{I_n\}$ is a *$b$-insensitivity sequence* if $I_n \gg b_n$ and

$$\sup_{x\geq I_n} \sup_{0\leq t\leq b_n} \left|\frac{F(x-t+\Delta)}{F(x+\Delta)} - 1\right| \to 0. \tag{3}$$

The next lemma shows that such a sequence can always be found if $F$ is a (locally) long-tailed distribution.

Lemma 2.2. *Let $\{b_n\}$ be a given sequence for which $b_n \to \infty$. We have $F \in \mathcal{L}_\Delta$ if and only if there exists a $b$-insensitivity sequence for $F$.*

*Truncation.* Motivated by the relationship between insensitivity and the class $\mathcal{L}_\Delta$, our next goal is to find a convenient way to think about the class of (locally) subexponential distributions $\mathcal{S}_\Delta$.

Given a sequence $\{b_n\}$, we call $\{h_n\}$ a *$b$-truncation sequence* for $F$ if

$$\lim_{K\to\infty} \limsup_{n\to\infty} \sup_{x\geq h_n} \frac{n\mathsf{P}\{S_2 \in x+\Delta, \xi_1, \xi_2 \in (-\infty, -Kb_n) \cup (h_n, \infty)\}}{F(x+\Delta)} = 0. \tag{4}$$

It is not hard to see that $n\overline{F}(h_n) = o(1)$ for any $b$-truncation sequence. We will see in Lemma 2.3(ii) below that a $b$-truncation sequence is often independent of $\{b_n\}$, in which case we simply say that $\{h_n\}$ is a *truncation sequence*. The reason for including the factor $n$ in the numerator is indicated in Section 2.2.

At first sight, this definition may raise several questions. The following lemma therefore provides motivation for the definition, and also shows that it can often be simplified. In Section 6, we present some tools to find good truncation sequences. For instance, as we show in Lemma 6.2, finding a truncation sequence is often not much different from checking a subexponentiality property; for this, standard techniques can be used.

Recall that a function $f$ is almost decreasing if $f(x) \asymp \sup_{y \geq x} f(y)$.

LEMMA 2.3. *Let $\{b_n\}$ be a natural-scale sequence.*

(i) *$F \in \mathcal{S}_\Delta$ if and only if $F \in \mathcal{L}_\Delta$ and there exists a $b$-truncation sequence for $F$.*

(ii) *If $x \mapsto F(x+\Delta)$ is almost decreasing, then $\{h_n\}$ can be chosen independently of $b$. Moreover, in that case, $\{h_n\}$ is a truncation sequence if and only if*

$$\lim_{n\to\infty} \sup_{x \geq h_n} \frac{n\mathsf{P}\{S_2 \in x+\Delta, \xi_1 > h_n, \xi_2 > h_n\}}{F(x+\Delta)} = 0.$$

*Small steps.* We next introduce the fourth and last sequence that plays a central role in this paper. For a given sequence $h \equiv \{h_n\}$, we call the sequence $\{J_n\}$ an *$h$-small-steps sequence* if

$$\lim_{n\to\infty} \sup_{x \geq J_n} \sup_{z \geq x} \frac{\mathsf{P}\{S_n \in z+\Delta, \xi_1 \leq h_n, \ldots, \xi_n \leq h_n\}}{nF(x+\Delta)} = 0. \tag{5}$$

Note that the inner supremum is always attained for $z = x$ if $T = \infty$. Moreover, in conjunction with the existence of a sequence for which (1) holds, (7) below shows that it is always possible to find a small-steps sequence for a subexponential distribution. Since it is often nontrivial to find a good $h$-small-steps sequence, Section 7 is entirely devoted to this problem.

*Main results.* The next theorem is our main result.

THEOREM 2.1. *Let $\{b_n\}$ be a natural-scale sequence, $\{I_n\}$ be a $b$-insensitivity sequence, $\{h_n\}$ be a $b$-truncation sequence and $\{J_n\}$ be an $h$-small-steps sequence. If $h_n = O(b_n)$ and $h_n \leq J_n$, we have*

$$\lim_{n\to\infty} \sup_{x \geq I_n + J_n} \left| \frac{\mathsf{P} S_n \in x+\Delta}{nF(x+\Delta)} - 1 \right| = 0.$$

The next subsection provides an outline of the proof of this theorem; the full proof is given in Sections 4 and 5. In all of the examples worked out in Sections 8 and 9, $I_n$ and $J_n$ are of different orders, and the boundary $I_n + J_n$ can be replaced by $\max(I_n, J_n)$. Our proof of the theorem, however, heavily relies on the additive structure given in the theorem.

In a variety of applications with $\mathsf{E}\{\xi\} = 0$, one wishes to conclude that $\mathsf{P}\{S_n \in na + \Delta\} \sim n\mathsf{P}\{\xi_1 \in na + \Delta\}$ for $a > 0$. As noted, for instance, by Doney [11] and S. Nagaev [34], it is thus of interest whether $na$ lies in the big-jump domain. Our next result shows that this can be concluded under minimal and readily verified conditions. The definition of $O$-regular variation is recalled in Appendix A; further details can be found in Chapter 2 of Bingham, Goldie and Teugels [3].

Corollary 2.1. *Assume that $\mathsf{E}\{\xi\} = 0$ and $\mathsf{E}\{|\xi|^\kappa\} < \infty$ for some $1 < \kappa \le 2$. Assume also that $x \mapsto F(x + \Delta)$ is almost decreasing and that $x \mapsto x^\kappa F(x + \Delta)$ either belongs to $\mathcal{S}d$ or is $O$-regularly varying. If furthermore*

$$\lim_{x\to\infty} \sup_{0 \le t \le x^{1/\kappa}} \left| \frac{F(x - t + \Delta)}{F(x + \Delta)} - 1 \right| = 0, \tag{6}$$

*then for any $a > 0$,*

$$\lim_{n\to\infty} \sup_{x \ge a} \left| \frac{\mathsf{P}\{S_n \in nx + \Delta\}}{n\mathsf{P}\{\xi_1 \in nx + \Delta\}} - 1 \right| = 0.$$

2.2. *Outline and idea of the proof of Theorem 2.1.* The first ingredient in the proof of Theorem 2.1 is the representation

$$\begin{aligned} &\mathsf{P}\{S_n \in x + \Delta\} \\ &\quad = \mathsf{P}\{S_n \in x + \Delta, B_1, \ldots, B_n\} + n\mathsf{P}\{S_n \in x + \Delta, \bar{B}_1, B_2, \ldots, B_n\} \\ &\qquad + \sum_{k=2}^{n} \binom{n}{k} \mathsf{P}\{S_n \in x + \Delta, \bar{B}_1, \ldots, \bar{B}_k, B_{k+1}, \ldots, B_n\}, \end{aligned} \tag{7}$$

where we set $B_i = \{\xi_i \le h_n\}$. To control the last term in this expression, we use a special exponential bound. Note that this bound is intrinsically different from Kesten's exponential bound (e.g., [14], Lemma 1.3.5), for which ramifications can be found in [42].

Lemma 2.4. *For $k \ge 2$, set*

$$\varepsilon_{\Delta,k}(n) \equiv \sup_{x \ge h_n} \frac{\mathsf{P}\{S_k \in x + \Delta, \xi_1 > h_n, \xi_2 > h_n, \ldots, \xi_k > h_n\}}{F(x + \Delta)}$$

*and*

$$\eta_{\Delta,k}(n, K) \equiv \sup_{x \ge h_n} \frac{\mathsf{P}\{S_k \in x + \Delta, \xi_2 < -Kb_n, \ldots, \xi_k < -Kb_n\}}{F(x + \Delta)}.$$

*Then we have* $\varepsilon_{\Delta,k}(n) \le \varepsilon_{\Delta,2}(n)^{k-1}$ *and* $\eta_{\Delta,k}(n,K) \le \eta_{\Delta,2}(n,K)^{k-1}$.

Our next result relies on this exponential bound, and shows that the sum in (7) is negligible when $\{h_n\}$ is a truncation sequence. In this argument, the factor $n$ in the numerator of (4) plays an essential role. The next lemma is inspired by Lemma 4 of Rozovskii [40].

LEMMA 2.5. *If* $F \in \mathcal{L}_\Delta$ *and* $n\varepsilon_{\Delta,2}(n) = o(1)$ *for some sequence* $\{h_n\}$, *then we have as* $n \to \infty$, *uniformly for* $x \in \mathbb{R}$,

$$\begin{aligned} &\mathsf{P}\{S_n \in x+\Delta\} \\ &\quad = \mathsf{P}\{S_n \in x+\Delta, \xi_1 \le h_n, \dots, \xi_n \le h_n\} \\ &\qquad + n\mathsf{P}\{S_n \in x+\Delta, \xi_1 > h_n, \xi_2 \le h_n, \dots, \xi_n \le h_n\}(1+o(1)). \end{aligned} \tag{8}$$

If $x$ is in the "small-steps domain," that is, if $x \ge J_n$, then the first term is small compared to $nF(x+\Delta)$. Therefore, proving Theorem 2.1 amounts to showing that the last term in (8) behaves like $nF(x+\Delta)$.

This is where insensitivity plays a crucial role. Intuitively, on the event $B_2, \dots, B_n$, $S_n - \xi_1$ stays on its natural scale: $|S_n - \xi_1| = O(b_n)$. Therefore, $S_n \in x+\Delta$ is roughly equivalent with $\xi_1 \in x \pm O(b_n) + \Delta$ on this event. In the "insensitive" domain ($x \ge I_n$), we know that $F(x \pm O(b_n)+\Delta) \approx F(x+\Delta)$, showing that the last term in (8) is approximately $nF(x+\Delta)$.

## 3. Proofs for Section 2.

In this section we prove all claims in Section 2 except for Theorem 2.1. Throughout many of the proofs, for convenience, we omit the mutual dependence of the four sequences. For instance, an insensitivity sequence should be understood as a $b$-insensitivity sequence for some given natural-scale sequence $\{b_n\}$.

Throughout this section, we use the notation of Lemma 2.4, and abbreviate $\varepsilon_{\Delta,k}(n)$ by $\varepsilon_k(n)$ if $T = \infty$. This is shortened further if $k = 2$; we then simply write $\varepsilon(n)$.

PROOF OF LEMMA 2.1. We derive a bound on $\mathsf{P}\{S_n > x, |\xi_1| \le cb_n, \dots, |\xi_n| \le cb_n\}$, which implies (by symmetry) the second estimate. A simple variant of the argument yields the first estimate.

Suppose that $\{S_n/b_n\}$ is tight. The first step in the proof is to show that

$$\lim_{K\to\infty} \liminf_{n\to\infty} \mathsf{P}\{S_n \in [-K^2 b_n, K^2 b_n], |\xi_1| \le Kb_n, \dots, |\xi_n| \le Kb_n\} = 1. \tag{9}$$

To see this, we observe that

$$\begin{aligned} &\mathsf{P}\{S_n \in [-K^2 b_n, K^2 b_n], |\xi_1| \le Kb_n, \dots, |\xi_n| \le Kb_n\} \\ &\quad \ge \mathsf{P}\{S_n \in [-K^2 b_n, K^2 b_n]\} - [1 - \mathsf{P}\{|\xi_1| \le Kb_n, \dots, |\xi_n| \le Kb_n\}]. \end{aligned}$$

By first letting $n$ tend to infinity and then $K$, we see that the first term tends to 1 by the tightness assumption, and the second term tends to zero by (2).

We next use a symmetrization argument. Let $S'_n$ be an independent copy of the random walk $S_n$, with step sizes $\xi'_1, \xi'_2, \ldots$. By (9), there exists a constant $K > 0$ such that $\mathsf{P}\{S'_n \le K^2 b_n, |\xi'_1| \le K b_n, \ldots, |\xi'_n| \le K b_n\} \ge 1/2$. On putting $\widetilde{S}_n = S_n - S'_n$ and $\widetilde{\xi}_i = \xi_i - \xi'_i$, we obtain

$$\begin{aligned}
&\mathsf{P}\{S_n > x, |\xi_1| \le cb_n, \ldots, |\xi_n| \le cb_n\} \\
&\quad \le 2\mathsf{P}\{S_n > x, |\xi_1| \le cb_n, \ldots, |\xi_n| \le cb_n, S'_n \le K^2 b_n, \\
&\qquad\qquad |\xi'_1| \le K b_n, \ldots, |\xi'_n| \le K b_n\} \\
&\quad \le 2\mathsf{P}\{\widetilde{S}_n > x - K^2 b_n, |\widetilde{\xi}_1| \le (c+K) b_n, \ldots, |\widetilde{\xi}_n| \le (c+K) b_n\}.
\end{aligned}$$

By the Chebyshev inequality, this is further bounded by

$$2\exp\left\{-sx + sK^2 b_n + n\log\int_{-(c+K)b_n}^{(c+K)b_n} e^{sz}\widetilde{F}(dz)\right\}$$

for all $s \ge 0$. Here, $\widetilde{F}$ denotes the distribution of $\xi_1 - \xi_2$. We use this inequality for $s = 1/(cb_n)$, implying that $sK^2 b_n$ is uniformly bounded in $n$ and $c \ge 1$. It remains to show that the same holds true for the last term in the exponent.

The key ingredient to bound this term is the assumption that $\{S_n/b_n\}$, and hence its symmetrized version $\{S'_n/b_n\}$, is tight. In the proof of the lemma in Section IX.7 of [17], Feller shows that there then exists some $c_0$ such that

$$A_0 \equiv \sup_n n\frac{\mathsf{E}\{\min(\widetilde{\xi}^2, (c_0 b_n)^2)\}}{b_n^2} < \infty.$$

It is convenient to also introduce $B_0 \equiv \sup_{y \le K+1}(e^y - 1 - y)/y^2$. In conjunction with the symmetry of $\widetilde{F}$, this immediately yields, for any $c \ge 1$,

$$\begin{aligned}
n\log\int_{-(c+K)b_n}^{(c+K)b_n} e^{sz}\widetilde{F}(dz) &\le n\int_{-(c+K)b_n}^{(c+K)b_n} e^{sz}\widetilde{F}(dz) - n \\
&\le n\int_{-(c+K)b_n}^{(c+K)b_n} [e^{sz} - 1 - sz]\widetilde{F}(dz) \\
&\le B_0 n\frac{\int_{-(c+K)b_n}^{(c+K)b_n} z^2\widetilde{F}(dz)}{c^2 b_n^2}.
\end{aligned}$$

Now, if $1 \le c < c_0 - K$, we bound this by $B_0 n b_n^{-2}\int_{-c_0 b_n}^{c_0 b_n} z^2\widetilde{F}(dz) \le A_0 B_0$. In the complementary case $c \ge c_0 - K$, we use the monotonicity of the function

$x \mapsto x^{-2}\mathsf{E}\{\min(\widetilde{\xi}^2, x^2)\}$ to see that

$$n\frac{\int_{-(c+K)b_n}^{(c+K)b_n} z^2 \widetilde{F}(dz)}{c^2 b_n^2} \le \frac{(c+K)^2}{c^2} n \frac{\mathsf{E}\{\min(\widetilde{\xi}^2, (c+K)^2 b_n^2)\}}{(c+K)^2 b_n^2}$$
$$\le (1+K)^2 n \frac{\mathsf{E}\{\min(\widetilde{\xi}^2, c_0^2 b_n^2)\}}{c_0^2 b_n^2},$$

which is bounded by $A_0(1+K)^2/c_0^2$. □

Proof of Lemma 2.2. Since $b_n \to \infty$, it is readily seen that $F \in \mathcal{L}_\Delta$ if $\{I_n\}$ is an insensitivity sequence. For the converse, we exploit the fact that $x \mapsto F(\log x + \Delta)$ is slowly varying. The uniform convergence theorem for slowly varying functions (see, e.g. Bingham, Goldie and Teugels [3], Theorem 1.2.1) implies that there exists some function $A$, increasing to $+\infty$, such that for $z \to \infty$,

$$\sup_{x \ge z} \sup_{0 \le y \le A(z)} \left| \frac{F(x-y+\Delta)}{F(x+\Delta)} - 1 \right| \to 0.$$

To complete the proof, it remains to choose $I_n = A^{-1}(b_n)$. □

Proof of Lemma 2.3. We shall first prove (ii), for which it is sufficient to show that $n\eta_{\Delta,2}(n,K)$ vanishes as first $n \to \infty$ and then $K \to \infty$. Observe that

$$n\mathsf{P}\{S_2 \in x+\Delta, \xi_2 < -Kb_n\} = n\int_{-\infty}^{-Kb_n} F(dy)F(x-y+\Delta)$$
$$\le nF(-Kb_n)\sup_{y \ge x} F(y+\Delta),$$

which is (up to multiplication by a finite constant) bounded by $nF(-Kb_n)F(x+\Delta)$ for large $x$ as $F(\cdot+\Delta)$ is almost decreasing. The claim therefore follows from (2).

Let us now prove (i). Let $F \in S_\Delta$. From $F \in \mathcal{L}_\Delta$ we deduce that we can find some function $h$ with $h(L) \le L/2$ and $h(L) \to \infty$ such that

$$\lim_{L\to\infty} \sup_{x \ge L} \sup_{y \in [-h(L), h(L)]} \left| \frac{F(x+y+\Delta)}{F(x+\Delta)} - 1 \right| = 0. \tag{10}$$

We start by showing that

$$\lim_{L\to\infty} \sup_{x \ge L} \frac{\mathsf{P}\{S_2 \in x+\Delta, \xi_1 > L, \xi_2 > L\}}{F(x+\Delta)}$$
$$= \lim_{L\to\infty} \sup_{x \ge 2L} \frac{\mathsf{P}\{S_2 \in x+\Delta, \xi_1 > L, \xi_2 > L\}}{F(x+\Delta)} = 0. \tag{11}$$

The first equality is only nontrivial if $T=\infty$, and can be deduced by considering $L \le x < 2L$ and $x \ge 2L$ separately. Next note that for $x \ge 2L$, since $h(2L) \le L$,

$$\begin{aligned}&\mathsf{P}\{S_2 \in x+\Delta, \xi_1 > L, \xi_2 > L\}\\ &\qquad \le \mathsf{P}\{S_2 \in x+\Delta\} - 2\mathsf{P}\{S_2 \in x+\Delta, \xi_2 \le h(2L)\}.\end{aligned} \tag{12}$$

We deduce (11) from the definitions of $h$ and $F \in \mathcal{S}_\Delta$.

In the global case $T=\infty$, (11) guarantees the existence of a truncation sequence for any $F \in \mathcal{S}_\Delta$ in view of part (ii) of the lemma. Slightly more work is required to prove this existence if $T<\infty$, relying on the bound

$$\begin{aligned}&\mathsf{P}\{S_2 \in x+\Delta, \xi_1, \xi_2 \in (-\infty, -Kb_n) \cup (h_n, \infty)\}\\ &\qquad \le \mathsf{P}\{S_2 \in x+\Delta, \xi_1, \xi_2 > h_n\}\\ &\qquad\quad + 2\mathsf{P}\{S_2 \in x+\Delta, \xi_1 < -Kb_n\}.\end{aligned} \tag{13}$$

As for the second term, we note that for any $x \gg b_n$

$$\begin{aligned}2\mathsf{P}\{S_2 \in x+\Delta, \xi_1 < -Kb_n\} &= 2\mathsf{P}\{S_2 \in x+\Delta, \xi_1 \le Kb_n\}\\ &\quad - 2\mathsf{P}\{S_2 \in x+\Delta, |\xi_1| \le Kb_n\}\\ &\le \mathsf{P}\{S_2 \in x\} - 2\mathsf{P}\{S_2 \in x+\Delta, |\xi_1| \le Kb_n\}.\end{aligned}$$

With (10), we readily find some $h_n \gg b_n$ such that $\mathsf{P}\{S_2 \in x+\Delta, |\xi_1| \le Kb_n\}$ is asymptotically equivalent to $G(Kb_n)F(x+\Delta)$ uniformly for $x \ge h_n$. We may assume without loss of generality that $n|\mathsf{P}S_2 \in x/F(x+\Delta) - 2| \to 0$ in this domain, so that by (2) the second term on the right-hand side of (13) is $o(1/n)F(x+\Delta)$ uniformly for $x \ge h_n$, as first $n \to \infty$ and then $K \to \infty$. In view of (11), we may also assume without loss of generality that $\mathsf{P}\{S_2 \in x+\Delta, \xi_1, \xi_2 > h_n\}$ is $o(1/n)F(x+\Delta)$ uniformly for $x \ge h_n$.

We have now shown that truncation sequences can be constructed if $F \in \mathcal{S}_\Delta$, and we proceed to the proof of the converse claim under the assumption $F \in \mathcal{L}$. Suppose that we are given some $\{h_n\}$ and $\{b_n\}$ such that (4) holds. For $x \ge 2h_n$, we have

$$0 \le \mathsf{P}\{S_2 \in x+\Delta\} - 2\mathsf{P}\{S_2 \in x+\Delta, \xi_1 \in [-Kb_n, h_n]\} \le \eta_{\Delta,2}(n,K)F(x+\Delta).$$

Again with (10), we readily find some $f_n \gg h_n$ such that $\mathsf{P}\{S_2 \in x+\Delta, \xi_1 \in [-Kb_n, h_n]\}$ is asymptotically equivalent to $F(x+\Delta)$ uniformly for $x \ge f_n$. Therefore $F \in \mathcal{S}_\Delta$. $\square$

Proof of Lemma 2.4. We only show that the first inequality holds; the second is simpler to derive and uses essentially the same idea.

Consider the global case $T=\infty$. We prove the inequality by induction. For $k=2$, the inequality is an equality. We now assume that the assertion holds for $k-1$ and we prove it for $k$. Recall that $B_j = \{\xi_j \le h_n\}$. First, for $x < kh_n$,

$$\begin{aligned}\mathsf{P}\{S_k > x, \bar{B}_1, \dots, \bar{B}_k\} &= \overline{F}(h_n)^k \le \varepsilon_{k-1}(n)\overline{F}((k-1)h_n)\overline{F}(h_n) \\ &\le \varepsilon_{k-1}(n)\varepsilon(n)\overline{F}(kh_n) \le \varepsilon(n)^{k-1}\overline{F}(x).\end{aligned}$$

Second, for $x \ge kh_n$,

$$\begin{aligned}&\mathsf{P}\{S_k > x, \bar{B}_1, \dots, \bar{B}_k\} \\ &\quad = \mathsf{P}\{\xi_k > x - h_n\}\mathsf{P}\{S_{k-1} > h_n, \bar{B}_1, \dots, \bar{B}_{k-1}\} \\ &\qquad + \int_{h_n}^{x-h_n} F(dz)\mathsf{P}\{S_{k-1} > x - z, \bar{B}_1, \dots, \bar{B}_{k-1}\} \\ &\quad \le \varepsilon_{k-1}(n)\left(\overline{F}(x-h_n)\overline{F}(h_n) + \int_{h_n}^{x-h_n} F(dz)\overline{F}(x-z)\right) \\ &\quad \le \varepsilon_{k-1}(n)\varepsilon(n)\overline{F}(x) \le \varepsilon(n)^{k-1}\overline{F}(x).\end{aligned}$$

This proves the assertion in the global case.

In the local case $T < \infty$, we again use induction. We may suppose that $h_n > T$. For $k = 2$, the claim is trivial. Assume now that it holds for $k-1$ and prove the inequality for $k$. First, for $x < kh_n - T$, it is clear that $\mathsf{P}\{\xi_1 > h_n, \xi_2 > h_n, \dots, \xi_k > h_n, S_k \in x + \Delta\} = 0$. Second, for $x \ge kh_n - T$,

$$\begin{aligned}&\mathsf{P}\{S_k \in x + \Delta, \xi_1 > h_n, \xi_2 > h_n, \dots, \xi_k > h_n\} \\ &\quad \le \int_{h_n}^{x-(k-1)h_n+T} F(dy) \\ &\qquad\qquad \times \mathsf{P}\{S_{k-1} \in x - y + \Delta, \xi_1 > h_n, \dots\} \\ &\quad \le \varepsilon_\Delta(n)^{k-2}\int_{h_n}^{x-(k-1)h_n+T} F(dy)F(x-y+\Delta) \\ &\quad \le \varepsilon_\Delta(n)^{k-2}\int_{h_n}^{x-h_n} F(dy)F(x-y+\Delta),\end{aligned}$$

where the latter inequality follows from the fact that $(k-1)h_n - T \ge h_n$ for $k > 2$. Now note that

$$\begin{aligned}&\int_{h_n}^{x-h_n} F(dy)F(x-y+\Delta) \\ &\quad \le \mathsf{P}\{S_2 \in x + \Delta, \xi_1 > h_n, \xi_2 > h_n\} \le \varepsilon_\Delta(n)F(x+\Delta),\end{aligned}$$

and the claim follows in the local case. □

We separately prove Lemma 2.5 in the global case and the local case.

Proof of Lemma 2.5: the global case. The assumption $F \in \mathcal{L}$ is not needed in the global case. For $k \geq 2$, we have

$$\begin{aligned}
&\mathsf{P}\{S_n > x, \bar{B}_1, \dots, \bar{B}_k, B_{k+1}, \dots, B_n\} \\
&\quad = \mathsf{P}\{\bar{B}_1, \dots, \bar{B}_k\}\mathsf{P}\{S_n - S_k > x - h_n, B_{k+1}, \dots, B_n\} \\
&\qquad + \mathsf{P}\{S_n > x, S_n - S_k \leq x - h_n, \bar{B}_1, \dots, \bar{B}_k, B_{k+1}, \dots, B_n\}.
\end{aligned}$$

We write $P_1$ and $P_2$ for the first and second summands respectively. Since $\overline{F}(h_n) \leq \varepsilon(n)$, the first term is estimated as follows:

$$P_1 \leq \varepsilon(n)^{k-1}\mathsf{P}\{S_n - S_k > x - h_n, \bar{B}_1, B_{k+1}, \dots, B_n\}.$$

Lemma 2.4 is used to bound the second term:

$$\begin{aligned}
P_2 &= \int_{-\infty}^{x-h_n} \mathsf{P}\{S_n - S_k \in dz, B_{k+1}, \dots, B_n\}\mathsf{P}\{S_k > x - z, \bar{B}_1, \dots, \bar{B}_k\} \\
&\leq \varepsilon(n)^{k-1} \int_{-\infty}^{x-h_n} \mathsf{P} S_n - S_k \in dz, B_{k+1}, \dots, B_n \overline{F}(x - z) \\
&= \varepsilon(n)^{k-1}\mathsf{P}\{\xi_1 + S_n - S_k > x, S_n - S_k \leq x - h_n, \bar{B}_1, B_{k+1}, \dots, B_n\}.
\end{aligned}$$

By combining these two estimates, we obtain that

$$\begin{aligned}
&\mathsf{P}\{S_n > x, \bar{B}_1, \dots, \bar{B}_k, B_{k+1} \dots, B_n\} \\
&\quad \leq \varepsilon(n)^{k-1}\mathsf{P}\{\xi_1 + S_n - S_k > x, \bar{B}_1, B_{k+1}, \dots, B_n\}.
\end{aligned}$$

Further,

$$\begin{aligned}
&\mathsf{P}\{S_n > x, \bar{B}_1, B_2, \dots, B_n\} \\
&\quad \geq \mathsf{P}\{S_n > x, \bar{B}_1, B_2, \dots, B_n, \xi_2 \geq 0, \dots, \xi_k \geq 0\} \\
&\quad \geq \mathsf{P}\{\xi_1 + S_n - S_k > x, \bar{B}_1, B_2, \dots, B_n, \xi_2 \geq 0, \dots, \xi_k \geq 0\} \\
&\quad = \mathsf{P}\{\xi_1 + S_n - S_k > x, \bar{B}_1, B_{k+1}, \dots, B_n\}\mathsf{P}\{0 \leq \xi_2 \leq h_n\}^{k-1}.
\end{aligned}$$

If $n$ is large enough, then $\mathsf{P}\{0 \leq \xi_2 \leq h_n\} \geq \mathsf{P}\{\xi_1 \geq 0\}/2 \equiv \beta$. Therefore, it follows from the above inequalities that

$$\begin{aligned}
&\mathsf{P}\{\xi_1 + S_n - S_k > x, \bar{B}_1, B_{k+1}, \dots, B_n\} \\
&\quad \leq \mathsf{P}\{S_n > x, \bar{B}_1, B_2, \dots, B_n\}\left(\frac{1}{\beta}\right)^{k-1}.
\end{aligned}$$

As a result, we have, for sufficiently large $n$,

$$\begin{aligned}
&\sum_{k=2}^{n} \binom{n}{k} \mathsf{P}\{S_n > x, \bar{B}_1, \dots, \bar{B}_k, B_{k+1}, \dots, B_n\} \\
&\quad \leq \mathsf{P}\{S_n > x, \bar{B}_1, B_2, \dots, B_n\} \sum_{k=2}^{n} \binom{n}{k} \left(\frac{\varepsilon(n)}{\beta}\right)^{k-1} \\
&\quad = o(n)\mathsf{P}\{S_n > x, \bar{B}_1, B_2, \dots, B_n\},
\end{aligned}$$

as desired. □

Proof of Lemma 2.5: the local case. We may assume that $h_n > T$ without loss of generality. The exponential bound of Lemma 2.4 shows that, for $k \geq 2$,

$$\begin{aligned}
&\mathsf{P}\{S_n \in x+\Delta, \bar{B}_1, \ldots, \bar{B}_k, B_{k+1}, \ldots, B_n\} \\
&\quad = \mathsf{P}\{S_n \in x+\Delta, S_n - S_k \leq x - h_n, \bar{B}_1, \ldots, \bar{B}_k, B_{k+1}, \ldots, B_n\} \\
&\quad = \int_{-\infty}^{x-h_n} \mathsf{P}\{S_n - S_k \in dz, B_{k+1}, \ldots, B_n\} \mathsf{P}\{S_k \in x - z + \Delta, \bar{B}_1, \ldots, \bar{B}_k\} \\
&\quad \leq \varepsilon_\Delta(n)^{k-1} \int_{-\infty}^{x-h_n} \mathsf{P} S_n - S_k \in dz, B_{k+1}, \ldots, B_n F(x - z + \Delta) \\
&\quad \leq \varepsilon_\Delta(n)^{k-1} \mathsf{P}\{\xi_1 + S_n - S_k \in x + \Delta, \bar{B}_1, B_{k+1}, \ldots, B_n\}.
\end{aligned}$$

Let $x_1 > 0$ be a constant such that $F(0, x_1] \equiv \beta > 0$. Then, for $n$ large enough so that $h_n > x_1$,

$$\begin{aligned}
&\mathsf{P}\{\xi_1 + S_n - S_k \in x + \Delta, \bar{B}_1, B_{k+1}, \ldots, B_n\} \\
&\quad = \beta^{1-k} \mathsf{P}\{\xi_1 + S_n - S_k \in x + \Delta, \bar{B}_1, B_{k+1}, \ldots, B_n, 0 < \xi_2, \ldots, \xi_k \leq x_1\} \\
&\quad \leq \beta^{1-k} \mathsf{P}\{S_n \in (x, x + (k-1)x_1 + T], \\
&\qquad\qquad \bar{B}_1, B_{k+1}, \ldots, B_n, 0 < \xi_2, \ldots, \xi_k \leq x_1\} \\
&\quad \leq \beta^{1-k} \mathsf{P}\{S_n \in (x, x + (k-1)x_1 + T], \bar{B}_1, B_2, \ldots, B_n\}.
\end{aligned}$$

Furthermore, we have

$$\begin{aligned}
&\mathsf{P}\{S_n \in (x, x + (k-1)x_1 + T], \bar{B}_1, B_2, \ldots, B_n\} \\
&\quad = \int_{-\infty}^{x-h_n+(k-1)x_1+T} \mathsf{P}\{\xi_1 > h_n, \xi_1 \in (x - y, x - y + (k-1)x_1 + T]\} \\
&\qquad\qquad \times \mathsf{P}\{S_n - \xi_1 \in dy, B_2, \ldots, B_n\}.
\end{aligned}$$

The condition $F \in \mathcal{L}_\Delta$ ensures that we can find some $x_0$ such that for any $x \geq x_0$, the inequality $F(x + T + \Delta) \leq 2F(x + \Delta)$ holds. Assuming without loss of generality that $x_1/T$ is an integer, this implies that for $y \leq x - h_n + (k-1)x_1 + T$ and $n$ large enough so that $h_n \geq x_0$,

$$\begin{aligned}
&\mathsf{P}\{\xi_1 > h_n, \xi_1 \in (x - y, x - y + (k-1)x_1 + T]\} \\
&\quad = \mathsf{P}\{\xi_1 \in (\max(h_n, x - y), x - y + (k-1)x_1 + T]\} \\
&\quad \leq \sum_{j=0}^{(k-1)x_1/T} F(\max(h_n, x - y) + jT + \Delta)
\end{aligned}$$

$$\leq \sum_{j=0}^{(k-1)x_1/T} 2^j F(\max(h_n, x-y)+\Delta)$$

$$\leq 2^{(k-1)x_1/T+1} F(\max(h_n, x-y)+\Delta).$$

Upon combining all inequalities that we have derived in the proof, we conclude that for large $n$, uniformly in $x \in \mathbb{R}$,

$$\mathsf{P}\{S_n \in x+\Delta, \bar{B}_1, \ldots, \bar{B}_k, B_{k+1}, \ldots, B_n\}$$
$$\leq 2\varepsilon_\Delta(n)^{k-1}\mathsf{P}\{S_n \in x+\Delta, \bar{B}_1, B_2, \ldots, B_n\}\left(\frac{2^{x_1/T}}{\beta}\right)^{k-1}.$$

The proof is completed in exactly the same way as for the global case. $\square$

Proof of Corollary 2.1. Let $a > 0$ be arbitrary, and note that it suffices to prove the claim for $a$ replaced by $2a$. By the Kolmogorov–Marcinkiewicz–Zygmund law of large numbers or the central limit theorem we can take $b_n = (na)^{1/\kappa}$. We readily check with (6) that $\{I_n \equiv an\}$ is an insensitivity sequence, and next show that $\{J_n \equiv an\}$ is a small-steps sequence. Observe that we may set $h_n = (na)^{1/\kappa}$ by Lemma 6.1 or Lemma 6.2 below. Therefore, we conclude with Lemma 2.1 that

$$\sup_{x\geq an}\sup_{z\geq x}\frac{\mathbf{P}(S_n \in z+\Delta, \xi_1 \leq (na)^{1/\kappa}, \ldots, \xi_n \leq (na)^{1/\kappa})}{F(x+\Delta)}$$
$$= O(1)\sup_{x\geq an}\frac{e^{-x^{1-1/\kappa}}}{F(x+\Delta)}.$$

Now we exploit the insensitivity condition (6) to prove that this upper bound vanishes. It implies that for any $\delta > 0$, there exists some $x_0 = x_0(\delta) > 0$ such that

$$\inf_{x\geq x_0}\frac{F(x+\Delta)}{F(x-x^{1/\kappa}+\Delta)} \geq 1-\delta.$$

In particular, $F(x+\Delta) \geq (1-\delta)^{x^{1-1/\kappa}}F(x/2+\Delta)$ for $x/2 \geq x_0$. Iterating, we obtain

$$\frac{F(x+\Delta)}{F(x_0+\Delta)} \geq (1-\delta)^{x^{1-1/\kappa}+(x/2)^{1-1/\kappa}+(x/4)^{1-1/\kappa}+\cdots} = e^{x^{1-1/\kappa}\ln(1-\delta)/(1-2^{-(1-1/\kappa)})}.$$

Since we may choose $\delta > 0$ small enough, we conclude that $e^{-x^{1-1/\kappa}} = o(F(x+\Delta))$ uniformly for $x \geq an$. It remains to apply Theorem 2.1. $\square$

**4. Proof of Theorem 2.1: the global case.** We separately prove the upper and lower bounds in Theorem 2.1, starting with the lower bound.

PROOF OF THEOREM 2.1: LOWER BOUND. For any $K > 0$ and $x \geq 0$, we have

$$\begin{aligned}
&\mathsf{P}\{S_n > x\} \\
&\quad \geq n\mathsf{P}\{S_n > x, \xi_1 > Kb_n, |\xi_2| \leq \sqrt{K}b_n, \ldots, |\xi_n| \leq \sqrt{K}b_n\} \\
&\quad \geq n\mathsf{P}\{\xi > x + Kb_n\}\mathsf{P}\{S_{n-1} > -Kb_n, |\xi_1| \leq \sqrt{K}b_n, \ldots, |\xi_{n-1}| \leq \sqrt{K}b_n\}.
\end{aligned}$$

Now let $\epsilon > 0$ be arbitrary, and fix some (large) $K$ such that

$$\begin{aligned}
&\liminf_{n\to\infty} \mathsf{P}\{S_{n-1} \in [-Kb_n, Kb_n], |\xi_1| \leq \sqrt{K}b_n, \ldots, |\xi_{n-1}| \leq \sqrt{K}b_n\} \\
&\quad \geq 1 - \epsilon/2,
\end{aligned} \tag{14}$$

which is possible by (9). Since $\{I_n\}$ is an insensitivity sequence, provided $n$ is large enough, we have $\overline{F}(x - b_n) \leq (1+\epsilon)^{1/K}\overline{F}(x)$ for any $x \geq I_n$. In particular, $\overline{F}(x + Kb_n) \geq (1+\epsilon)^{-1}\overline{F}(x)$ for $x \geq I_n$. Conclude that for any $x \geq I_n$,

$$\frac{\mathsf{P}\{S_n > x\}}{n\mathsf{P}\{\xi > x\}} \geq (1+\epsilon)^{-1}\mathsf{P}\{S_{n-1} > -Kb_n, |\xi_1| \leq \sqrt{K}b_n, \ldots, |\xi_{n-1}| \leq \sqrt{K}b_n\},$$

which must exceed $(1+\epsilon)^{-1}(1-\epsilon)$ for large enough $n$. $\square$

PROOF OF THEOREM 2.1: UPPER BOUND. Since $\{J_n\}$ is a small-steps sequence, it suffices to focus on the second term on the right-hand side of (8).

Fix some (large) $K$, and suppose throughout that $x \geq I_n + J_n$. Recall that $B_i = \{\xi_i \leq h_n\}$. Since $I_n \gg b_n$ and $h_n = O(b_n)$, we must have $x - J_n \geq h_n$ for large $n$. We may therefore write

$$\begin{aligned}
&\mathsf{P}\{S_n > x, \bar{B}_1, B_2, \ldots, B_n\} \\
&\qquad = \int_{h_n}^{x-J_n} + \int_{x-J_n}^{\infty} F(du)\mathsf{P}\{S_n - \xi_1 > x - u, B_2, \ldots, B_n\}.
\end{aligned} \tag{15}$$

For $u$ in the first integration interval, we clearly have $x - u \geq J_n$, so that by construction of $\{J_n\}$ and $\{h_n\}$, for large $n$,

$$\begin{aligned}
&\int_{h_n}^{x-J_n} F(du)P(S_{n-1} > x - u, B_1, \ldots, B_{n-1}) \\
&\qquad \leq e^{-K}n\int_{h_n}^{x-J_n} F(du)\overline{F}(x-u) \leq e^{-K}n\int_{h_n}^{x-h_n} F(du)\overline{F}(x-u) \\
&\qquad \leq e^{-K}n\mathsf{P}\{S_2 > x, \xi_1 > h_n, \xi_2 > h_n\} \leq e^{-K}\overline{F}(x),
\end{aligned}$$

where we also used the assumption $J_n \geq h_n$.

In order to handle the second integral in (15), we rely on the following fact. As $\{I_n\}$ is an insensitivity sequence, we have for large $n$,

$$\sup_{u \geq I_n} \frac{\overline{F}(u)}{\overline{F}(u+b_n)} \leq e^{1/K^2}. \tag{16}$$

We next distinguish between two cases: $J_n \leq K b_n$ and $J_n > K b_n$. In the first case, since $x - J_n \geq I_n$, (16) can be applied iteratively to see that

$$\overline{F}(x - J_n) \leq e^{J_n/(K^2 b_n)} \overline{F}(x) \leq e^{1/K} \overline{F}(x). \tag{17}$$

Now note that the second integral in (15) is majorized by $\mathsf{P}\{\xi > x - J_n\}$ and hence by $e^{1/K}\overline{F}(x)$.

Slightly more work is needed if $J_n > K b_n$. First write the last integral in (15) as $\int_{x-J_n}^{x-Kb_n} + \int_{x-Kb_n}^{\infty}$. Since $x - Kb_n > x - J_n \geq I_n$, the argument of the preceding paragraph shows that $\mathsf{P}\{\xi > x - Kb_n\} \leq e^{1/K}\overline{F}(x)$. This must also be an upper bound for the integral $\int_{x-Kb_n}^{\infty}$, so it remains to investigate the integral $\int_{x-J_n}^{x-Kb_n}$, which is bounded from above by

$$\begin{aligned} &\mathsf{P}\{\xi > x - J_n\}\mathsf{P}\{S_{n-1} > J_n, B_1, \dots, B_{n-1}\} \\ &\qquad + \int_{Kb_n}^{J_n} \mathsf{P}\{S_{n-1} \in dy, B_1, \dots, B_{n-1}\}\overline{F}(x-y). \end{aligned}$$

First, using $h_n = O(b_n)$, select some $c < \infty$ such that $h_n \leq c b_n$. Without loss of generality, we may suppose that $K^2 > c$. Using the first inequality in (17) and Lemma 2.1, we see that the first term is bounded by $O(1)e^{J_n/(K^2 b_n) - J_n/(c b_n)} \times \overline{F}(x) = o(1)\overline{F}(x)$ as $n \to \infty$. As $x - J_n \geq I_n$, the second term is bounded by

$$\begin{aligned} &\sum_{k=K}^{\lfloor J_n/b_n \rfloor} \mathsf{P}\{S_{n-1}/b_n \in (k, k+1], \xi_1 \leq h_n, \dots, \xi_{n-1} \leq h_n\}\overline{F}(x - kb_n) \\ &\quad \leq \sum_{k=K}^{\lfloor J_n/b_n \rfloor} \mathsf{P}\{S_{n-1} > kb_n, \xi_1 \leq cb_n, \dots, \xi_{n-1} \leq cb_n\}\overline{F}(x - kb_n) \\ &\quad \leq C \sum_{k=K}^{\lfloor J_n/b_n \rfloor} e^{-k/c} e^{k/K^2} \overline{F}(x) \leq C \frac{e^{-K/c + 1/K}}{1 - e^{-1/c + 1/K^2}} \overline{F}(x), \end{aligned}$$

where we have used (16) and (the first inequality of) Lemma 2.1. Since $K$ is arbitrary, this proves the upper bound. $\square$

## 5. Proof of Theorem 2.1: the local case.

We use the following notation throughout this section: set $C_i^K \equiv \{-\sqrt{K} b_n \leq \xi_i \leq h_n\}$ and $D_i^K \equiv \{\xi_i <$ $-\sqrt{K}b_n\}$ for any $K>0$. Recall that $B_i=\{\xi_i\le h_n\}$. As in Section 4, we start with the lower bound.

PROOF OF THEOREM 2.1: LOWER BOUND. The proof is similar to its global analogue, again using (14) and insensitivity. First fix some $\epsilon>0$, then choose $K$ (fixed) as in the "global" proof. For later use, by (2) we may assume without loss of generality that $K$ satisfies $\sup_n n\overline{G}(Kb_n)<\epsilon$ and that $e^{-1/K}\ge 1-\epsilon$.

Repeated application of "insensitivity" shows that for any $y\ge 0$, provided $n$ is large,

$$\inf_{x\ge I_n}\frac{F(x+y+\Delta)}{F(x+\Delta)}\ge \exp\left\{-\frac{y}{K^2b_n}\right\},\qquad \sup_{x\ge I_n}\frac{F(x+y+\Delta)}{F(x+\Delta)}\le \exp\left\{\frac{y}{K^2b_n}\right\}.$$

We next distinguish between the cases $J_n\ge Kb_n$ and $J_n<Kb_n$. In the first case, since we consider $x\ge I_n+J_n$, we have $x-Kb_n\ge I_n$ for large $n$, so that

$$\begin{aligned}
&\mathsf{P}\{S_n\in x+\Delta\}\\
&\quad\ge n\int_{-Kb_n}^{Kb_n}\mathsf{P}\{S_{n-1}\in dy,|\xi_1|\le\sqrt{K}b_n,\dots,|\xi_{n-1}|\le\sqrt{K}b_n\}F(x-y+\Delta)\\
&\quad\ge ne^{-1/K}\mathsf{P}\{S_{n-1}\in[-Kb_n,Kb_n],|\xi_1|\le\sqrt{K}b_n,\dots,|\xi_{n-1}|\le\sqrt{K}b_n\}\\
&\qquad\times F(x+\Delta)\\
&\quad\ge ne^{-1/K}(1-\epsilon)F(x+\Delta),
\end{aligned}$$

where the second inequality uses the above insensitivity relations (distinguish between positive and negative $y$). Since $e^{-1/K}\ge 1-\epsilon$, this proves the claim if $J_n\ge Kb_n$.

We next suppose that $J_n<Kb_n$. Observe that then, for $x\ge I_n+J_n$,

$$\inf_{-J_n\le y\le 0}\frac{F(x+y+\Delta)}{F(x+\Delta)}\ge\exp\left\{-\frac{J_n}{K^2b_n}\right\}\ge e^{-1/K}.$$

Since $h_n=O(b_n)$ and $I_n\gg b_n$, the events $C_1^K$ and $\{\xi_1>x-J_n\}$ are disjoint for $x\ge I_n+J_n$, so that with the preceding display,

$$\begin{aligned}
&\mathsf{P}\{S_n\in x+\Delta\}\\
&\quad\ge n\int_{-Kb_n}^{J_n}\mathsf{P}\{S_{n-1}\in dy,C_1^K,\dots,C_{n-1}^K\}F(x-y+\Delta)\\
&\quad\ge ne^{-1/K}F(x+\Delta)\mathsf{P}\{S_{n-1}\in[-Kb_n,J_n],C_1^K,\dots,C_{n-1}^K\}\\
&\quad\ge n(1-\epsilon)F(x+\Delta)\mathsf{P}\{S_{n-1}\in[-Kb_n,J_n],C_1^K,\dots,C_{n-1}^K\}.
\end{aligned}$$

We need two auxiliary observations before proceeding. First, by construction of $K$, we have

$$\begin{aligned}&\mathsf{P}\{S_{n-1}<-Kb_n,C_1^K,\dots,C_{n-1}^K\}\\&\qquad\le\mathsf{P}\{|S_{n-1}|>Kb_n,|\xi_1|\le\sqrt{K}b_n,\dots,|\xi_{n-1}|\le\sqrt{K}b_n\}\le\epsilon.\end{aligned}$$

Furthermore, by definition of $J_n$, we have for large $n$,

$$\begin{aligned}&\mathsf{P}\{S_{n-1}>J_n,C_1^K,\dots,C_{n-1}^K\}\\&\qquad\le\mathsf{P}\{S_{n-1}>J_n,\xi_1\le h_n,\dots,\xi_{n-1}\le h_n\}\\&\qquad=\sum_{k=0}^{\infty}\mathsf{P}\{S_{n-1}\in J_n+kT+\Delta,\xi_1\le h_n,\dots,\xi_{n-1}\le h_n\}\\&\qquad\le\epsilon n\sum_{k=0}^{\infty}F(J_n+kT+\Delta)=\epsilon n\overline{F}(J_n)\le\epsilon n\overline{F}(h_n)\le\epsilon,\end{aligned}$$

since $n\overline{F}(h_n)=o(1)$.

The inequalities in the preceding two displays show that

$$\mathsf{P}\{S_{n-1}\in[-Kb_n,J_n],C_1^K,\dots,C_{n-1}^K\}\ge\mathsf{P}\{C_1^K\}^{n-1}-2\epsilon,$$

and by construction of $K$ we may infer that $\mathsf{P}\{C_1^K\}\ge 1-\overline{F}(h_n)-\overline{G}(Kb_n)\ge 1-2\epsilon/n$, so that $\mathsf{P}\{C_1^K\}^{n-1}$ must exceed $e^{-3\epsilon}$ if $n$ is large. □

The proof of the upper bound is split into two lemmas, Lemma 5.1 and Lemma 5.2. First note that by Lemma 2.5 and the definition of $J_n$, it suffices to show that

$$\limsup_{n\to\infty}\sup_{x\ge I_n+J_n}\frac{\mathsf{P}\{S_n\in x+\Delta,\xi_1>h_n,\xi_2\le h_n,\dots,\xi_n\le h_n\}}{F(x+\Delta)}\le 1.$$

We prove this by truncation from below. The numerator in the preceding display can be rewritten as

$$\begin{aligned}&\mathsf{P}\{S_n\in x+\Delta,\bar{B}_1,C_2^K,\dots,C_n^K\}\\&\qquad+\sum_{k=2}^{n}\binom{n-1}{k-1}\mathsf{P}\{S_n\in x+\Delta,\bar{B}_1,D_2^K,\dots,D_k^K,C_{k+1}^K,\dots,C_n^K\}.\end{aligned}\tag{18}$$

The first probability in this expression is taken care of by the next lemma.

Lemma 5.1. *Under the assumptions of Theorem 2.1, we have*

$$\limsup_{K\to\infty}\limsup_{n\to\infty}\sup_{x\ge I_n+J_n}\frac{\mathsf{P}\{S_n\in x+\Delta,\bar{B}_1,C_2^K,\dots,C_n^K\}}{F(x+\Delta)}\le 1.$$

Proof. This is similar to the "global" proof of Theorem 2.1, but some new arguments are needed. We follow the lines of the proof given in Section 4.

Fix some (large) $K>1$. Suppose that $n$ is large enough such that

$$\sup_{x\ge I_n}\frac{F(x+b_n+\Delta)}{F(x+\Delta)}\le e^{1/K^2}. \tag{19}$$

In order to bound the probability

$$\begin{aligned}&\mathsf{P}\{S_n\in x+\Delta, h_n<\xi_1\le x-\min(J_n,Kb_n),C_2^K,\dots,C_n^K\}\\&\quad\le \mathsf{P}\{S_n\in x+\Delta, h_n<\xi_1\le x-\min(J_n,Kb_n),B_2,\dots,B_n\},\end{aligned}$$

exactly the same arguments work as for the global case.

Moreover, after distinguishing between $J_n>Kb_n$ and $J_n\le Kb_n$, it is not hard to see with (19) that for $x\ge I_n+J_n$ and $n$ large,

$$\begin{aligned}&\mathsf{P}\{S_n\in x+\Delta, x-\min(J_n,Kb_n)<\xi_1\le x+Kb_n,C_2^K,\dots,C_n^K\}\\&\quad=\int_{-Kb_n}^{\min(J_n,Kb_n)+T}\mathsf{P}\{S_{n-1}\in dy,C_1^K,\dots,C_{n-1}^K\}F(x-y+\Delta)\\&\quad\le e^{1/K}\mathsf{P}\{S_{n-1}\in[-Kb_n,\min(J_n,Kb_n)+T],C_1^K,\dots,C_{n-1}^K\}F(x+\Delta),\end{aligned}$$

which is majorized by $e^{1/K}F(x+\Delta)$.

It remains to investigate the regime $\xi_1>x+Kb_n$. Since $h_n=O(b_n)$, we may assume without loss of generality that $h_n\le\sqrt{K}b_n$. Exploiting the insensitivity inequality (19) and the second inequality of Lemma 2.1, we see that for $x\ge I_n$ and $n$ large enough,

$$\begin{aligned}&\mathsf{P}\{S_n\in x+\Delta,\xi_1>x+Kb_n,C_2^K,\dots,C_n^K\}\\&\quad\le\int_{-\infty}^{T-Kb_n}\mathsf{P}\{S_{n-1}\in dy,|\xi_1|\le\sqrt{K}b_n,\dots,|\xi_{n-1}|\le\sqrt{K}b_n\}\\&\qquad\qquad\times F(x-y+\Delta)\\&\quad\le e^{1/K^2}\sum_{k=K-1}^{\infty}\mathsf{P}\{|S_{n-1}|>kb_n,|\xi_1|\le\sqrt{K}b_n,\dots,|\xi_{n-1}|\le\sqrt{K}b_n\}\\&\qquad\qquad\times F(x+kb_n+\Delta)\\&\quad\le Ce^{1/K^2}\sum_{k=K-1}^{\infty}e^{-k/\sqrt{K}}e^{k/K^2}F(x+\Delta)\\&\quad=Ce^{1/K^2}\frac{e^{-(K-1)/\sqrt{K}+(K-1)/K^2}}{1-e^{-1/\sqrt{K}+1/K^2}}F(x+\Delta).\end{aligned}$$

It is not hard to see (e.g., with l'Hôpital's rule) that the prefactor can be made arbitrarily small. □

The next lemma deals with the sum over $k$ in (18). Together with Lemma 5.1, it completes the proof of Theorem 2.1 in the local case.

LEMMA 5.2. *Under the assumptions of Theorem 2.1,*

$$\limsup_{n\to\infty} \sup_{x\geq I_n+J_n} \frac{\sum_{k=2}^n \binom{n-1}{k-1} \mathsf{P}\{S_n\in x+\Delta, \bar{B}_1, D_2^K,\dots,D_k^K, C_{k+1}^K,\dots,C_n^K\}}{F(x+\Delta)}$$

*converges to zero as $k\to\infty$.*

PROOF. The $k$th term in the sum can be written as

$$\begin{aligned}
&\mathsf{P}\{S_n\in x+\Delta, \bar{B}_1, D_2^K,\dots,D_k^K, C_{k+1}^K,\dots,C_n^K\}\\
&\quad = \mathsf{P}\{S_n\in x+\Delta, \bar{B}_1, D_2^K,\dots,D_k^K, C_{k+1}^K,\dots,C_n^K, S_n-S_k\leq x-h_n\} \qquad (20)\\
&\qquad + \mathsf{P}\{S_n\in x+\Delta, \bar{B}_1, D_2^K,\dots,D_k^K, C_{k+1}^K,\dots,C_n^K, S_n-S_k> x-h_n\}.
\end{aligned}$$

As for the first term, we know that by definition of $\eta_{\Delta,k}$,

$$\begin{aligned}
&\mathsf{P}\{S_n\in x+\Delta, S_n-S_k\leq x-h_n, \bar{B}_1, D_2^K,\dots,D_k^K, C_{k+1}^K,\dots,C_n^K\}\\
&\quad = \int_{-\infty}^{x-h_n} \mathsf{P}\{S_n-S_k\in dy, C_{k+1}^K,\dots,C_n^K\}\\
&\qquad\qquad \times \mathsf{P}\{S_k\in x-y+\Delta, \bar{B}_1, D_2^K,\dots,D_k^K\}\\
&\quad \leq \eta_{\Delta,k}(n,\sqrt{K})\mathsf{P}\{\xi_1+S_n-S_k\in x+\Delta, \bar{B}_1, C_{k+1}^K,\dots,C_n^K\}.
\end{aligned}$$

The arguments of the proof of Lemma 2.5 in the local case can be repeated to see that there exists some $\gamma>0$ independent of $K$, $n$ and $x$, such that for any $x$,

$$\begin{aligned}
&\mathsf{P}\{\xi_1+S_n-S_k\in x+\Delta, \bar{B}_1, C_{k+1}^K,\dots,C_n^K\}\\
&\qquad \leq 2\gamma^{k-1}\mathsf{P}\{S_n\in x+\Delta, \bar{B}_1, C_2^K,\dots,C_n^K\}.
\end{aligned}$$

As $n\to\infty$ and then $K\to\infty$, the probability on the right-hand side is bounded by $F(x+\Delta)$ in view of Lemma 5.1. We use the assumption on $\eta_{\Delta,2}(n,K)$ to study the prefactor: with Lemma 2.4 and some elementary estimates, we obtain

$$\lim_{K\to\infty}\limsup_{n\to\infty}\sum_{k=2}^n \binom{n-1}{k-1}\gamma^{k-1}\eta_{\Delta,k}(n,\sqrt{K})=0.$$

We now proceed to the second term on the right-hand side of (20):

$$\begin{aligned}
&\mathsf{P}\{S_n \in x+\Delta, S_n - S_k > x - h_n, \bar{B}_1, D_2^K, \ldots, D_k^K, C_{k+1}^K, \ldots, C_n^K\} \\
&\quad \le \int_{-\infty}^{h_n+T} \mathsf{P}\{S_k \in dy, \bar{B}_1, D_2^K, \ldots, D_k^K\} \\
&\qquad\qquad \times \mathsf{P}\{S_n - S_k \in x - y + \Delta, C_{k+1}^K, \ldots, C_n^K\} \\
&\quad \le \mathsf{P}\{\bar{B}_1, D_2^K, \ldots, D_k^K\} \sup_{z > x - h_n - T} \mathsf{P}\{S_n - S_k \in z + \Delta, C_{k+1}^K, \ldots, C_n^K\}.
\end{aligned}$$

Since $\{b_n\}$ and $\{h_n\}$ are natural-scale and truncation sequences, respectively, the first probability is readily seen to be $o(n^{-k})$ as first $n \to \infty$ and then $K \to \infty$.

In order to investigate the supremum in the preceding display, we choose $x_0 > 0$ such that $F(x_0 + \Delta) \equiv \beta > 0$. Without loss of generality, we may assume that $h_n > x_0$. Then we have

$$\begin{aligned}
&\mathsf{P}\{S_n - S_k \in z + \Delta, C_{k+1}^K, \ldots, C_n^K\} \\
&\quad = \beta^{-k} \mathsf{P}\{S_n - S_k \in z + \Delta, C_{k+1}^K, \ldots, C_n^K, \xi_1 \in x_0 + \Delta, \ldots, \xi_k \in x_0 + \Delta\} \\
&\quad \le \beta^{-k} \mathsf{P}\{S_n \in z + kx_0 + (k+1)\Delta, \\
&\qquad\qquad C_{k+1}^K, \ldots, C_n^K, \xi_1 \in x_0 + \Delta, \ldots, \xi_k \in x_0 + \Delta\} \\
&\quad \le \beta^{-k} \sum_{j=0}^{k} \mathsf{P}\{S_n \in z + kx_0 + jT + \Delta, C_1^K, \ldots, C_n^K\} \\
&\quad \le 2k\beta^{-k} \sup_{u > z} \mathsf{P}\{S_n \in u + \Delta, C_1^K, \ldots, C_n^K\},
\end{aligned}$$

showing that

$$\begin{aligned}
&\sup_{z > x - h_n - T} \mathsf{P}\{S_n - S_k \in z + \Delta, C_{k+1}^K, \ldots, C_n^K\} \\
&\qquad \le 2k\beta^{-k} \sup_{z > x - h_n - T} \mathsf{P}\{S_n \in z + \Delta, C_1^K, \ldots, C_n^K\}.
\end{aligned}$$

This implies that, uniformly for $x \ge I_n + J_n$, as $n \to \infty$ and then $K \to \infty$,

$$\begin{aligned}
&\sum_{k=2}^{n} \binom{n-1}{k-1} \mathsf{P}\{S_n \in x + \Delta, \bar{B}_1, D_2, \ldots, D_k, C_{k+1}^K, \ldots, C_n^K, S_n - S_k > x - h_n\} \\
&\quad = \sum_{k=2}^{n} \binom{n-1}{k-1} o(n^{-k}) k\beta^{-k} \sup_{z > x - h_n - T} \mathsf{P}\{S_n \in z + \Delta, C_1^K, \ldots, C_n^K\} \\
&\quad = o(1/n) \sup_{z > x - h_n - T} \mathsf{P}\{S_n \in z + \Delta, C_1^K, \ldots, C_n^K\}
\end{aligned}$$

$$\le o(1/n) \sup_{z>x-h_n-T} \mathsf{P}\{S_n \in z+\Delta, \xi_1 \le h_n, \dots, \xi_n \le h_n\}$$
$$\le o(1) F(x-h_n-T+\Delta),$$

where we have used the definition of the small-steps sequence $\{J_n\}$, in conjunction with the assumptions that $h_n = O(b_n)$ and $I_n \gg b_n$.

Since $J_n \ge h_n$, we clearly have $x - h_n \ge I_n$ in the regime $x \ge I_n + J_n$. Therefore, insensitivity shows that $F(x-h_n-T+\Delta) = O(1)F(x+\Delta)$, and the claim follows. $\square$

**6. On truncation sequences.** It is typically nontrivial to choose good truncation and small-steps sequences. Therefore, we devote the next two sections to present some techniques which are useful for selecting $\{h_n\}$ and $\{J_n\}$. The present section focuses on truncation sequences $\{h_n\}$. We give two tools for selecting truncation sequences.

We first investigate how to choose a truncation sequence in the presence of $O$-regular variation (see Appendix A).

LEMMA 6.1. *If $x \mapsto F(x+\Delta)$ is almost decreasing and $O$-regularly varying, then $\{h_n\}$ is a truncation sequence if $n\overline{F}(h_n) = o(1)$.*

PROOF. Let us first suppose that $T = \infty$. Using Lemma 2.3(ii), the claim is proved once we have shown that $\varepsilon_{\Delta,2}(n) = o(1/n)$ if $n\overline{F}(h_n) = o(1)$. To this end, we write

$$\mathsf{P}\{S_2 > x, \xi_1 > h_n, \xi_2 > h_n\} \le 2\mathsf{P}\{\xi_1 > x/2, \xi_2 > h_n\} = 2\overline{F}(h_n)\overline{F}(x/2),$$

and note that for $x \ge h_n$, $\overline{F}(x/2) = O(\overline{F}(x))$ as a result of the assumption that $\overline{F}$ is $O$-regularly varying.

For the local case, it suffices to prove that $n\varepsilon_{\Delta,2} = o(1)$ if $n\overline{F}(h_n) = o(1)$. Since the mapping $x \mapsto F(x+\Delta)$ is $O$-regularly varying, the uniform convergence theorem for this class (Theorem 2.0.8 in [3]) implies that $\sup_{y\in[1/2,1]} F(xy+\Delta) \le CF(x+\Delta)$ for some constant $C < \infty$ (for large enough $x$). Therefore, if $n$ is large, we have for $x \ge h_n$,

$$\mathsf{P}\{S_2 \in x+\Delta, \xi_1 > h_n, \xi_2 > h_n\}$$
$$\le 2\mathsf{P}\{S_2 \in x+\Delta, h_n < \xi_1 \le x/2+T, \xi_2 > x/2\}$$
$$\le 2\int_{h_n}^{x/2+T} F(dy)F(x-y+\Delta),$$

which is bounded by $2C\overline{F}(h_n)F(x+\Delta)$; the claim follows. $\square$

The next lemma is our second tool for selecting truncation sequences. For the definition of $\mathcal{S}d$, we refer to Appendix B.

LEMMA 6.2. *Let $x \mapsto F(x+\Delta)$ be almost decreasing, and suppose that $x \mapsto x^r F(x+\Delta)$ belongs to $\mathcal{S}d$ for some $r > 0$. Then any $\{h_n\}$ with $\limsup_{n\to\infty} n h_n^{-r} < \infty$ is a truncation sequence.*

PROOF. Set $H(x) \equiv x^r F(x+\Delta)$, and first consider $T = \infty$. It follows from $F \in \mathcal{L}$ that for large $n$

$$\begin{aligned}\int_{h_n}^{x/2} \overline{F}(x-y)F(dy) &\le \sum_{i=\lfloor h_n\rfloor}^{\lceil x/2\rceil} \overline{F}(x-i)F(i,i+1] \\ &\le \sum_{i=\lfloor h_n\rfloor}^{\lceil x/2\rceil} \overline{F}(x-i)\overline{F}(i) \\ &\le 2\sum_{i=\lfloor h_n\rfloor}^{\lceil x/2\rceil} \int_i^{i+1} \overline{F}(x-y)\overline{F}(y)\,dy \\ &\le 2\int_{h_n}^{x/2} \overline{F}(x-y)\overline{F}(y)\,dy.\end{aligned}$$

By Lemma 2.3(ii) and the above arguments, we obtain

$$\begin{aligned}\varepsilon_{\Delta,2}(n) &= \sup_{x\ge 2h_n} \frac{\overline{F}(x/2)^2 + 2\int_{h_n}^{x/2}\overline{F}(x-y)F(dy)}{\overline{F}(x)} \\ &\le 2\sup_{x\ge 2h_n} \frac{\overline{F}(x/2)^2 + 2\int_{h_n}^{x/2}\overline{F}(x-y)\overline{F}(y)\,dy}{\overline{F}(x)} \\ &\le \frac{2^{r+1}}{h_n^r}\sup_{x\ge 2h_n}\left(\frac{H(x/2)^2 + \int_{h_n}^{x/2} H(x-y)H(y)\,dy}{H(x)}\right).\end{aligned}$$

We now exploit the assumption that $H \in \mathcal{S}d$. First observe that $\int_{x/2-T}^{x/2} H(y)H(x-y)\,dy = o(H(x))$ if $H \in \mathcal{S}d$, implying $H(x/2)^2 = o(H(x))$ in conjunction with $H \in \mathcal{L}$. We deduce that for any $M > 0$,

$$\varepsilon_{\Delta,2}(n) \le o(h_n^{-r}) + O(h_n^{-r})\sup_{x\ge 2h_n} \frac{\int_M^{x/2} H(y)H(x-y)\,dy}{H(x)},$$

so that $\varepsilon_{\Delta,2}(n) = o(h_n^{-r})$.

Let us now turn to the case $T < \infty$. Note that by Lemma 2.3(ii), we exploit the long-tailedness of $x \mapsto F(x+\Delta)$ to obtain, for large $n$,

$$\varepsilon_{\Delta,2}(n) \le \sup_{x\ge 2h_n - T} \frac{2\int_{h_n}^{(x+T)/2} F(x-y+\Delta)F(dy)}{F(x+\Delta)}$$

$$\le 4 \sup_{x\ge 2h_n} \frac{\int_{h_n}^{x/2} F(x-y+\Delta)F(dy)}{F(x+\Delta)}.$$

An elementary approximation argument, again relying on the long-tailedness assumption, shows that uniformly for $x \ge 2h_n$,

$$\int_{h_n}^{x/2} F(x-y+\Delta)F(dy) \sim \frac{1}{T}\int_{h_n}^{x/2} F(y+\Delta)F(x-y+\Delta)\,dy.$$

The rest of the proof parallels the global case. □

**7. On small-steps sequences.** In this section, we investigate techniques that are often useful for selecting small-steps sequences $\{J_n\}$. That is, we derive bounds on $\mathsf{P}\{S_n \in x+\Delta, \xi_1 \le h_n, \dots, \xi_n \le h_n\}$ under a variety of assumptions.

We first need some more notation. Write $\varphi_n = \mathsf{E}\{e^{\xi/h_n}; \xi \le h_n\}$, and let $\{\xi_i^{(n)}\}_{i=1}^\infty$ be a sequence of "twisted" (or "tilted") i.i.d. random variables with distribution function

$$\mathsf{P}\{\xi^{(n)} \le y\} = \frac{\mathsf{E}\{e^{\xi/h_n}; \xi \le h_n, \xi \le y\}}{\varphi_n}.$$

We also put $S_k^{(n)} = \xi_1^{(n)} + \cdots + \xi_k^{(n)}$; note that $\{S_k^{(n)}\}$ is a random walk for any $n$.

Next we introduce a sequence $\{a_n\}$ which plays an important role in the theory of domains of (partial) attraction. First define $Q(x) \equiv x^{-2}\mu_2(x) + \overline{G}(x)$. It is not hard to see that $Q$ is continuous, ultimately decreasing and that $Q(x) \to 0$ as $x \to \infty$. Therefore, the solution to the equation $Q(x) = n^{-1}$, which we call $a_n$, is well defined and unique for large $n$.

LEMMA 7.1. *We have the following exponential bounds.*

(i) *If $\mathsf{E}\{\xi\} = 0$ and $\mathsf{E}\{\xi^2\} = 1$, then for any $\epsilon > 0$ there exists some $n_0$ such that for any $n \ge n_0$ and any $x \ge 0$,*

$$\mathsf{P}\{S_n \in x+\Delta, \xi_1 \le h_n, \dots, \xi_n \le h_n\} \le \exp\left\{-\frac{x}{h_n} + \left(\frac{1}{2}+\epsilon\right)\frac{n}{h_n^2}\right\}\mathsf{P}\{S_n^{(n)} \in x+\Delta\}.$$

(ii) *If $h_n \ge a_n$ and $n|\mu_1(a_n)| = O(a_n)$, then there exists some $C < \infty$ such that for any $n \ge 1$ and any $x \ge 0$,*

$$\mathsf{P}\{S_n \in x+\Delta, \xi_1 \le h_n, \dots, \xi_n \le h_n\} \le C\exp\left\{-\frac{x}{h_n}\right\}\mathsf{P}\{S_n^{(n)} \in x+\Delta\}.$$

(iii) *If* $\mathsf{E}\{\xi\} = 0$ *and* $x \mapsto F(-x)$ *is regularly varying at infinity with index* $-\alpha$ *for some* $\alpha \in (1,2)$*, then for any* $\epsilon > 0$ *there exists some* $n_0$ *such that for any* $n \ge n_0$ *and any* $x \ge 0$,

$$\begin{aligned}&\mathsf{P}\{S_n \in x+\Delta, \xi_1 \le h_n, \ldots, \xi_n \le h_n\}\\ &\quad\le \exp\biggl\{-\frac{x}{h_n} + \frac{n}{h_n^2}\int_0^{h_n} u^2 F(du) + (1+\epsilon)\frac{\Gamma(2-\alpha)}{\alpha-1} nF(-h_n)\biggr\}\\ &\qquad\times \mathsf{P}\{S_n^{(n)} \in x+\Delta\}.\end{aligned}$$

(iv) *If* $x \mapsto F(-x)$ *is regularly varying at infinity with index* $-\alpha$ *for some* $\alpha \in (0,1)$*, then for any* $\epsilon > 0$ *there exists some* $n_0$ *such that for any* $n \ge n_0$ *and any* $x \ge 0$,

$$\begin{aligned}&\mathsf{P}\{S_n \in x+\Delta, \xi_1 \le h_n, \ldots, \xi_n \le h_n\}\\ &\quad\le \exp\biggl\{-\frac{x}{h_n} + \frac{n}{h_n}\int_0^{h_n} uF(du)\\ &\qquad\qquad + \frac{n}{h_n^2}\int_0^{h_n} u^2F(du) - (1-\epsilon)\Gamma(1-\alpha)nF(-h_n)\biggr\}\\ &\qquad\times \mathsf{P}\{S_n^{(n)} \in x+\Delta\}.\end{aligned}$$

PROOF. We need to investigate $n\log\varphi_n$ under the four sets of assumptions of the lemma, since

$$\begin{aligned}\mathsf{P}\{S_n \in x+\Delta, \xi_1 \le h_n, \ldots, \xi_n \le h_n\} &= \varphi_n^n \mathsf{E}\{e^{-S_n^{(n)}/h_n}; S_n^{(n)} \in x+\Delta\}\\ &\le e^{-x/h_n + n\log\varphi_n}\mathsf{P}\{S_n^{(n)} \in x+\Delta\}.\end{aligned}$$

We start with (i). Since $e^y \le 1 + y + y^2/2 + |y|^3$ for $y \le 1$, some elementary bounds in the spirit of the proof of Lemma 2.1 show that

$$n\log\varphi_n \le n\int_{-h_n}^{h_n}[e^{z/h_n} - 1]F(dz) \le \frac{n\mu_1(h_n)}{h_n} + \frac{n\mu_2(h_n)}{2h_n^2} + \frac{n\overline{\mu}_3(h_n)}{h_n^3},$$

where $\overline{\mu}_3(h_n) = \int_{-h_n}^{h_n}|z|^3F(dz)$. It follows from $\mathsf{E}\{\xi^2\} = 1$ that $\overline{\mu}_3(h_n) = o(h_n)$. Indeed, if $\mathsf{E}\{\xi^2\} < \infty$, then $\mathsf{E}\{\xi^2 f(|\xi|)\} < \infty$ for some function $f(x) \uparrow \infty, x/f(x) \uparrow \infty$, so that

$$\overline{\mu}_3(h_n) = \int_{-h_n}^{h_n}|z|^3F(dz) \le \frac{h_n}{f(h_n)}\int_{-h_n}^{h_n} z^2 f(z)F(dz) = O(1)\frac{h_n}{f(h_n)} = o(h_n).$$

One similarly gets $\mu_1(h_n) = o(1/h_n)$, relying on $\mathsf{E}\{\xi\} = 0$. This proves the first claim.

For (ii), we use similar arguments and the inequality $e^y - 1 \le y + y^2$ for $y \le 1$. From $h_n \ge a_n$ it follows that

$$\frac{n\mu_1(h_n)}{h_n} + \frac{n\mu_2(h_n)}{h_n^2} \le \frac{n|\mu_1(a_n)| + n\int_{a_n}^{h_n} yF(dy)}{h_n} + nQ(h_n)$$
$$\le \frac{n|\mu_1(a_n)|}{a_n} + n\overline{F}(a_n) + nQ(h_n).$$

The first term is bounded by assumption and the other two are both bounded by $nQ(a_n) = 1$.

To prove the third claim, we use $\mathsf{E}\{\xi\} = 0$ to write

$$n\log\varphi_n \le n\int_{-\infty}^{h_n} (e^{u/h_n} - 1 - u/h_n)F(du)$$
$$= n\biggl(\int_{-\infty}^{0} + \int_{0}^{h_n}\biggr)(e^{u/h_n} - 1 - u/h_n)F(du).$$

After integrating the first integral by parts twice, we see that

$$\int_{-\infty}^{0} (e^{u/h_n} - 1 - u/h_n)F(du) = h_n^{-2}\int_{-\infty}^{0} e^{u/h_n}\biggl(\int_{-\infty}^{u} F(t)\,dt\biggr)\,du.$$

By Karamata's theorem, $u \mapsto \int_{-\infty}^{-u} F(t)\,dt$ is regularly varying at infinity with index $-\alpha + 1$. We can thus apply a Tauberian theorem (e.g., [3], Theorem 1.7.1) to obtain for $n \to \infty$,

$$h_n^{-2}\int_{-\infty}^{0} e^{u/h_n}\biggl(\int_{-\infty}^{u} F(t)\,dt\biggr)\,du \sim h_n^{-1}\Gamma(2-\alpha)\int_{-\infty}^{-h_n} F(t)\,dt$$
$$\sim \frac{\Gamma(2-\alpha)}{\alpha - 1}F(-h_n).$$

We finish the proof of the third claim by observing that

$$\int_{0}^{h_n} (e^{u/h_n} - 1 - u/h_n)F(du) \le h_n^{-2}\int_{0}^{h_n} u^2 F(du).$$

Part (iv) is proved similarly, relying on the estimate

$$n\log\varphi_n \le n\biggl(\int_{-\infty}^{0} + \int_{0}^{h_n}\biggr)(e^{u/h_n} - 1)F(du).$$

After integrating the first integral by parts and applying a Tauberian theorem, we obtain

$$n\int_{-\infty}^{0} (e^{u/h_n} - 1)F(du) = -nh_n^{-1}\int_{-\infty}^{0} e^{u/h_n}F(u)\,du \sim -\Gamma(1-\alpha)nF(-h_n).$$

The integral over $[0, h_n]$ can be bounded using the inequality $e^y - 1 \le y + y^2$ for $y \le 1$. $\square$

In order to apply the estimates of the preceding lemma, we need to study $\mathsf{P}\{S_n^{(n)} \in x+\Delta\}$. If $T=\infty$, it is generally sufficient to bound this by 1, but in the local case we need to study our "truncated" and "twisted" random walk $\{S_k^{(n)}\}$ in more detail. Therefore, we next give a *concentration* result in the spirit of Gnedenko's local limit theorem. However, we do not restrict ourselves to distributions belonging to a domain of attraction. Instead, we work within the more general framework of Griffin, Jain and Pruitt [18] and Hall [19]. Our proof is highly inspired by these works, as well as by ideas of Esseen [15], Feller [16] and Petrov [36].

We need the following condition introduced by Feller [16]:

$$\limsup_{x\to\infty} \frac{x^2\overline{G}(x)}{\mu_2(x)} < \infty, \tag{21}$$

which also facilitates the analysis in [18, 19]. This condition is discussed in more detail in Section 9.1.

Proposition 7.1. *Suppose that we have either:*

1. $\mathsf{E}\{\xi^2\}<\infty$ *and* $\mathsf{E}\{\xi\}=0$, *or*
2. $\mathsf{E}\{\xi^2\}=\infty$ *and (21) holds.*

*Let $T<\infty$. There exist finite constants $C, C'$ such that, for all large $n$,*

$$\sup_{x\in\mathbb{R}} \mathsf{P}\{S_n^{(n)} \in x+\Delta\} \le \frac{C}{h_n} + \frac{C'}{a_n}.$$

Proof. Throughout, $C$ and $C'$ denote strictly positive, finite constants that may vary from line to line.

Let $\xi_s^{(n)}$ denote the symmetrized version of $\xi^{(n)}$, that is, $\xi_s^{(n)} = \xi_1^{(n)} - \xi_2^{(n)}$, where the $\xi_i^{(n)}$ are independent. For any $\epsilon>0$, we have the Esseen bound (see Petrov [36], Lemma 1.16 for a ramification)

$$\sup_{x\in\mathbb{R}} \mathsf{P}\{S_n^{(n)} \in x+\Delta\} \le C\epsilon^{-1}\int_{-\epsilon}^{\epsilon} |\mathsf{E}\{e^{it\xi^{(n)}}\}|^n\,dt.$$

Since $x \le \exp[-(1-x^2)/2]$ for $0\le x\le 1$ and $|\mathsf{E}\{e^{it\xi^{(n)}}\}|^2 = \mathsf{E}\{\cos t\xi_s^{(n)}\}$, this is further bounded by

$$\begin{aligned}
&C\epsilon^{-1}a_n^{-1}\int_{-\epsilon a_n}^{\epsilon a_n} |\mathsf{E}\{e^{it\xi^{(n)}/a_n}\}|^n\,dt \\
&\quad\le C\epsilon^{-1}a_n^{-1}\int_0^{\epsilon a_n} \exp[-(n/2)\mathsf{E}\{1-\cos(t\xi_s^{(n)}/a_n)\}]\,dt \\
&\quad\le C\epsilon^{-1}h_n^{-1} + C\epsilon^{-1}a_n^{-1}\int_{a_n/h_n}^{\epsilon a_n} \exp[-(n/2)\mathsf{E}\{1-\cos(t\xi_s^{(n)}/a_n)\}]\,dt.
\end{aligned}$$

Now note that for $h_n^{-1} \le t \le \epsilon$, provided $\epsilon$ is chosen small enough,

$$\begin{aligned}
\mathsf{E}\{1-\cos(t\xi_s^{(n)})\} &\ge Ct^2\mathsf{E}\{(\xi_s^{(n)})^2; |\xi_s^{(n)}| \le t^{-1}\} \\
&\ge C\varphi_n^{-2}t^2 \int_{\substack{x,y\le h_n \\ |x-y|\le t^{-1}}} (x-y)^2 e^{(x+y)/h_n} F(dx)F(dy) \\
&\ge C\varphi_n^{-2}t^2 \int_{|x|\le t^{-1}/2, |y|\le t^{-1}/2} (x-y)^2 e^{(x+y)/h_n} F(dx)F(dy) \\
&\ge Ct^2 e^{-t^{-1}/h_n}[\mu_2(t^{-1}/2) - \mu_1(t^{-1}/2)^2] \\
&\ge Ct^2\mu_2(t^{-1}/2) - Ct^2\mu_1(t^{-1}/2)^2.
\end{aligned}$$

If $\lim_{x\to\infty}\mu_2(x) < \infty$ and $\lim_{x\to\infty}\mu_1(x) = 0$, then it is clear that we can select $\epsilon$ so that, uniformly for $t \le \epsilon$,

$$\mu_2(t^{-1}/2) - \mu_1(t^{-1}/2)^2 \ge \mu_2(t^{-1}/2)/2.$$

The same can be done if $\mu_2(x) \to \infty$. Indeed, let $a > 0$ satisfy $\overline{G}(a) \le 1/8$. Application of the Cauchy–Schwarz inequality yields for $t < 1/(2a)$,

$$\begin{aligned}
\mu_1(t^{-1}/2)^2 &= (\mu_1(t^{-1}/2) - \mu_1(a) + \mu_1(a))^2 \\
&\le 2(\mu_1(t^{-1}/2) - \mu_1(a))^2 + 2\mu_1(a)^2 \\
&\le 2\mu_2(t^{-1}/2)\overline{G}(a) + 2\mu_1(a)^2 \le \mu_2(t^{-1}/2)/4 + 2\mu_1(a),
\end{aligned}$$

and the assumption $\mu_2(x) \to \infty$ shows that we can select $\epsilon$ small enough so that this is dominated by $\mu_2(t^{-1}/2)/2$ for $t \le \epsilon$.

Having seen that $\mathsf{E}\{1-\cos(t\xi_s^{(n)})\} \ge Ct^2\mu_2(t^{-1}/2)$, we next investigate the truncated second moment. To this end, we use (21), which always holds if $\mathsf{E}\{\xi^2\} < \infty$, to see that there exists some $C'$ such that $t^2\mu_2(t^{-1}/2)/2 \ge C'Q(t^{-1}/2)$.

We conclude that there exist some $\epsilon, C, C' \in (0,\infty)$ such that

$$\sup_{x\in\mathbb{R}} \mathsf{P}\{S_n^{(n)} \in x + \Delta\} \le C\epsilon^{-1}h_n^{-1} + C\epsilon^{-1}a_n^{-1}\int_{2a_n/h_n}^{2\epsilon a_n} \exp[-C'nQ(a_nt^{-1})]\,dt.$$

In order to bound the integral, we use the following result due to Hall [19]. Under (21), there exists some $k \ge 1$ such that for large enough $n$,

$$\int_k^{2\epsilon a_n} \exp[-C'nQ(a_nt^{-1})]\,dt \le C.$$

If $2a_n/h_n \ge k$, this immediately proves the claim. In the complementary case, we bound the integral over $[2a_n/h_n, k]$ simply by $k$. $\square$

**8. Examples with finite variance.** After showing heuristically how $J_n$ can be chosen, this section applies our main result (Theorem 2.1) to random walks with step-size distributions satisfying $\mathsf{E}\{\xi\}=0$ and $\mathsf{E}\{\xi^2\}=1$. Clearly, $\{S_n/\sqrt{n}\}$ is then tight and thus one can always take $b_n=\sqrt{n}$ as a natural-scale sequence.

Our goals are to show that our theory recovers many known large-deviation results, and that it fills gaps in the literature allowing new examples to be worked out. In fact, finding big-jump domains with our theory often essentially amounts to verifying whether the underlying step-size distribution is subexponential.

8.1. *A heuristic for choosing $J_n$.* Before showing how $J_n$ can typically be guessed in the finite-variance case, we state an auxiliary lemma of which the proof contains the main idea for the heuristic. Observe that the function $g$ in the lemma tends to infinity as a consequence of the finite-variance assumption.

LEMMA 8.1. *Consider $F$ for which $\mathsf{E}\{\xi\}=0$ and $\mathsf{E}\{\xi^2\}=1$. Let $g$ satisfy $-\log[x^2F(x+\Delta)]\le g(x)$ for large $x$ and suppose that $g(x)/x$ is eventually nonincreasing.*

*Let a sequence $\{J_n\}$ be given.*

(i) *If $T=\infty$, suppose that*

$$\limsup_{n\to\infty}\frac{g(J_n)}{J_n^2/n}<\frac{1}{2}. \tag{22}$$

(ii) *If $T<\infty$, suppose that*

$$\limsup_{n\to\infty}\frac{g(J_n)}{J_n^2/n+\log n}<\frac{1}{2}.$$

*If $\{h_n=n/J_n\}$ is a truncation sequence, then $\{J_n\}$ is a corresponding small-steps sequence.*

PROOF. Let $\epsilon>0$ be given. First consider the case $T=\infty$. By Lemma 7.1(i), we have to show that the given $h_n$ and $J_n$ satisfy

$$\sup_{x\ge J_n}\left(-\frac{x}{h_n}+\left(\frac{1}{2}+\epsilon\right)\frac{n}{h_n^2}-\log\overline{F}(x)-\log n\right)\to-\infty. \tag{23}$$

Next observe that $J_n\gg\sqrt{n}$, for otherwise $g(J_n)$ would be bounded; this is impossible in view of the assumption on $J_n$. Therefore, not only $g(x)/x$ is nondecreasing for $x\ge J_n$, but the same holds true for $\log[x^2/n]/x$. This yields, on substituting $h_n=n/J_n$,

$$\sup_{x\ge J_n}\left(-\frac{x}{h_n}-\log\overline{F}(x)-\log n\right)\le\sup_{x\ge J_n}x\left(-\frac{J_n}{n}+\frac{g(J_n)}{J_n}+\frac{\log[J_n^2/n]}{J_n}\right),$$

and the supremum is attained at $J_n$ since the expression between brackets is negative as a result of our assumption on $J_n$. Conclude that the left-hand side of (23) does not exceed

$$-\frac{1-\epsilon}{2}\frac{J_n^2}{n} + g(J_n) - \log\frac{J_n^2}{n},$$

which tends to $-\infty$ if $\epsilon$ is chosen appropriately.

The local case $T < \infty$ is similar. By Proposition 7.1 and Lemma 7.1(i), it now suffices to show

$$\sup_{x \geq J_n}\left(-\frac{x}{h_n} + \left(\frac{1}{2}+\epsilon\right)\frac{n}{h_n^2} - \log F(x+\Delta) - \log n - \log h_n\right) \to -\infty.$$

With the above arguments and the identity $2\log h_n = \log n - \log(J_n^2/n)$, it follows that the expression on the left-hand side is bounded by

$$-\frac{1-\epsilon}{2}\left[\frac{J_n^2}{n} + \log n\right] + g(J_n) - \frac{1}{2}\log\frac{J_n^2}{n},$$

and the statement thus follows from the assumption on $J_n$ as before. □

The idea of the above proof allows to heuristically find the *best possible* small-steps sequence in the finite-variance case. Let us work this out for $T = \infty$. Use (23) to observe that $J_n$ is necessarily larger than or equal to

$$\left(\frac{1}{2}+\epsilon\right)\frac{n}{h_n} - h_n \log n - h_n \log \overline{F}(J_n).$$

Set $\epsilon = 0$ for simplicity, and then minimize the right-hand side with respect to $h_n$. We find that the minimizing value (i.e., the best possible truncation sequence) is

$$h_n = \sqrt{\frac{n}{-2\log[n\overline{F}(J_n)]}}.$$

Since $h_n = n/J_n$ according to the above lemma, this suggests that the following asymptotic relation holds for the best small-steps sequence:

$$J_n \sim \sqrt{-2n\log[n\overline{F}(J_n)]}. \tag{24}$$

We stress that a number of technicalities need to be resolved before concluding that any $J_n$ satisfying this relation constitutes a small-steps sequence; the heuristic should be treated with care. In fact, one typically needs that $J_n$ is *slightly* bigger than suggested by (24). Still, we encourage the reader to compare the heuristic big-jump domain with the big-jump domain that we find for the examples in the remainder of this section.

8.2. *$O$-regularly varying tails.* In this subsection, it is our aim to recover A. Nagaev's classical boundary for regularly varying tails from Theorem 2.1. In fact, we work in the more general setting of $O$-regular variation.

PROPOSITION 8.1. *Suppose that $\mathsf{E}\{\xi\} = 0$ and $\mathsf{E}\{\xi^2\} = 1$. Moreover, let $x \mapsto F(x+\Delta)$ be $O$-regularly varying with upper Matuszewska index $\alpha_F$ and lower Matuszewska index $\beta_F$.*

1. *If $T = \infty$, suppose that $\alpha_F < -2$, and let $t > -\beta_F - 2$.*
2. *If $T < \infty$, suppose that $\alpha_F < -3$, and let $t > -\beta_F - 3$.*

*The sequence $\{h_n \equiv \sqrt{n/(t \log n)}\}$ is a truncation sequence. Moreover, for this choice of $h$, $\{J_n \equiv \sqrt{tn \log n}\}$ is an $h$-small-steps sequence.*

PROOF. We first show that $\{h_n\}$ is a truncation sequence, for which we use the third part of Lemma 2.3. In the global case, Theorem 2.2.7 in [3] implies that for any $\epsilon > 0$, we have $\overline{F}(x) \le x^{\alpha_F + \epsilon}$ for large $x$. By choosing $\epsilon$ small enough, we get $n\overline{F}(h_n) = o(1)$ since $\alpha_F < -2$. For the local case, we first need to apply Theorem 2.6.3(a) in [3] and then the preceding argument; this yields that for any $\epsilon > 0$, $\overline{F}(x) \le x^{1+\alpha_F+\epsilon}$ provided $x$ is large. Then we use $\alpha_F < -3$ to choose $\epsilon$ appropriately.

Our next aim is to show that $\{J_n\}$ is a small-steps sequence. We only do this for $T = \infty$; the complementary case is similar. Fix some $\epsilon > 0$ to be determined later. Again by Theorem 2.2.7 in [3], we know that $-\log[x^2\overline{F}(x)]$ is dominated by $(-2 - \beta_F + \epsilon) \log x$, which is eventually nonincreasing on division by $x$. Application of Lemma 8.1 shows that it suffices to choose an $\epsilon > 0$ satisfying

$$\limsup_{n\to\infty} \frac{(-2 - \beta_F + \epsilon) \log J_n}{J_n^2/n} < \frac{1}{2},$$

and it is readily seen that this can be done for $J_n$ given in the proposition. $\square$

With the preceding proposition at hand, we next derive the Nagaev boundary from Theorem 2.1. Indeed, as soon as an insensitivity sequence $\{I_n\}$ is determined, we can conclude that $\mathsf{P}\{S_n \in x + \Delta\} \sim nF(x + \Delta)$ uniformly for $x \ge I_n + J_n$, where the sequence $\{J_n\}$ is given in Proposition 8.1. Since $J_n$ depends on some $t$ which can be chosen appropriately, the above asymptotic equivalence holds uniformly for $x \ge J_n$ if $J_n \gg I_n$.

A class of distributions for which we can immediately conclude that $J_n \gg I_n$ is constituted by the requirement that $x \mapsto F(x + \Delta)$ is intermediate regularly varying (see Appendix A). Then any $I_n \gg b_n$ can be chosen as an insensitivity sequence; see Corollary 2.2I in [8].

The next theorem is due to A. Nagaev in the global case with regularly varying $\overline{F}$; see [14], Theorem 8.6.2 or Ng et al. [35]. In the local regularly varying case, it goes at least back to Pinelis [37].

THEOREM 8.1. *Let the assumptions of Proposition 8.1 hold, and suppose that $x \mapsto F(x+\Delta)$ is intermediate regularly varying at infinity.*

*With $t$ chosen as in Proposition 8.1, we have $\mathsf{P}\{S_n \in x+\Delta\} \sim nF(x+\Delta)$ uniformly for $x \geq \sqrt{tn\log n}$.*

8.3. *Logarithmic hazard function.* In this subsection, we consider step-size distributions with

$$F(x+\Delta) = p(x)e^{-c\log^\beta x},$$

where $\beta > 1, c > 0$ and $p$ is $O$-regularly varying with $p \in \mathcal{L}$. Note that log-normal distributions as well as Benktander Type I step-size distributions fit into this framework. Lemma B.1 with $R(x) = z(x) = c\log^\beta x$ shows that $x \mapsto F(x+\Delta)$ belongs to the class $\mathcal{S}d$ of subexponential densities.

We first select a small-steps sequence.

PROPOSITION 8.2. *Suppose that $\mathsf{E}\{\xi\} = 0$ and $\mathsf{E}\{\xi^2\} = 1$, and consider the above setup. Let $t > 2^{1-\beta}c$.*

*The sequence $\{h_n \equiv \sqrt{n/(t\log^\beta n)}\}$ is a truncation sequence, and $\{J_n \equiv \sqrt{tn\log^\beta n}\}$ is a corresponding small-steps sequence.*

PROOF. We only consider the global case, since the same arguments are used in the local case.

The family of distributions we consider is closed under multiplication by a polynomial. Moreover, $x \mapsto F(x+\Delta)$ is almost decreasing. To see this, write $F(x+\Delta) = p(x)x^\eta x^{-\eta}e^{-c\log^\beta x}$ and choose $\eta \in \mathbb{R}$ so that $p(x)x^\eta$ is almost decreasing; this can be done since the upper Matuszewska index of $p$ is finite. Membership of $\mathcal{S}d$ in conjunction with Lemma 6.2 shows that $\{h_n\}$ is a truncation sequence.

To show that $\{J_n\}$ is a corresponding small-steps sequence, we note that $p(x) \leq x^{c'}$ for some $c' \in \mathbb{R}$ provided $x$ is large [3], Theorem 2.2.7. We next use Lemma 8.1 with $g(x) = c'\log x + c\log^\beta x$. □

Theorem 2.1 yields the big-jump domain as soon as an insensitivity sequence is selected. This is readily done if $p$ is intermediate regularly varying, and we now work out this special case. First note that

$$\frac{F(x-\sqrt{n}+\Delta)}{F(x+\Delta)} = \frac{p(x-\sqrt{n})}{p(x)}\exp(c[\log^\beta x - \log^\beta(x-\sqrt{n})]).$$

Next observe that, by the uniform convergence theorem for regularly varying functions [3], Theorem 1.5.2, $x \gg \sqrt{n}$,

$$\log^\beta x - \log^\beta(x - \sqrt{n}) \le \beta\sqrt{n} \sup_{x-\sqrt{n}\le y\le x} y^{-1}\log^{\beta-1} y$$
$$\sim \beta x^{-1}\sqrt{n}\log^{\beta-1} x,$$

and a matching lower bound is derived similarly. This shows that, although the ratio of the $p$-functions converges uniformly to 1 in the domain $x \gg \sqrt{n}$, the analogous domain for the $\log^\beta$-functions is smaller. We conclude that any $I_n$ with $\sqrt{n}\log^{\beta-1} I_n = o(I_n)$ is an insensitivity sequence; in particular we may choose any $I_n$ satisfying $I_n \gg \sqrt{n\log^{2\beta-2} n}$.

We have thus proved the following theorem, which is new in the local case. As noted in [30], the "global" part (ii) can be deduced from Lemma 2A in Rozovskii [40]. The first part should be compared to Corollary 1 of [40], where a partial result is obtained.

THEOREM 8.2. *Let the assumptions of Proposition 8.2 hold, and suppose that $p$ is intermediate regularly varying at infinity.*

*With $t$ chosen as in Proposition 8.2, we have $\mathsf{P}\{S_n \in x + \Delta\} \sim nF(x + \Delta)$:*

(i) *uniformly for $x \ge \sqrt{tn\log^\beta n}$ if $1 < \beta < 2$, and*

(ii) *uniformly for $x \ge x_n$ for any $x_n \gg \sqrt{n\log^{2\beta-2} n}$ if $\beta \ge 2$.*

8.4. *Regularly varying hazard function.* In this subsection, we consider step-size distributions with

$$F(x + \Delta) = p(x)e^{-R(x)},$$

where $R$ is differentiable. We suppose that $p$ is $O$-regularly varying with $p \in \mathcal{L}$, and that $R'$ is regularly varying with index $\beta - 1$ for some $\beta \in (0, 1)$. In particular, by Karamata's theorem, $R$ is regularly varying with index $\beta$. Note that Weibull as well as Benktander Type II step-size distributions fit into this framework. Moreover, Lemma B.1 with $z(x) = x^\alpha$ for some $\alpha \in (\beta, 1)$ implies that $x \mapsto F(x + \Delta)$ belongs to $\mathcal{S}d$.

PROPOSITION 8.3. *Suppose that $\mathsf{E}\{\xi\} = 0$ and $\mathsf{E}\{\xi^2\} = 1$, and consider the above setup.*

*For any $\epsilon > 0$, the sequence $\{h_n \equiv n^{(1-\beta-\epsilon)/(2-\beta)}\}$ is a truncation sequence, and $\{J_n \equiv n^{(1+\epsilon)/(2-\beta)}\}$ is a corresponding small-steps sequence.*

PROOF. Along the lines of the proof of Proposition 8.2. In Lemma 8.1, we use $g(x) = x^{\beta+\epsilon^2}$. $\square$

In the above proposition, we have not given the best possible small-steps sequence, as any insensitivity sequence is asymptotically larger than $J_n$. To see this when $p$ is intermediate regularly varying, note that for $x \gg \sqrt{n}$

$$\frac{F(x-\sqrt{n}+\Delta)}{F(x+\Delta)} = \frac{p(x-\sqrt{n})}{p(x)} e^{R(x)-R(x-\sqrt{n})} \leq (1+o(1)) e^{\sqrt{n}\sup_{x-\sqrt{n}\leq y\leq x} R'(y)}.$$

Since $R'$ is regularly varying, we have $\sup_{x-\sqrt{n}\leq y\leq x} R'(y) \sim R'(x)$ if $x \gg \sqrt{n}$. A lower bound is proved along the same lines. The observation $R'(x) \asymp x^{-1}R(x)$ allows to show that $I_n \gg J_n$, and the next theorem follows on applying Theorem 2.1.

THEOREM 8.3. *Let the assumptions of Proposition 8.3 hold, and suppose that $p$ is intermediate regularly varying at infinity.*

*For any $\{x_n\}$ with $x_n/R(x_n) \gg \sqrt{n}$, we have $\mathsf{P}\{S_n \in x+\Delta\} \sim nF(x+\Delta)$ uniformly for $x \geq x_n$.*

8.5. *"Light" subexponential tails.* In this subsection, we consider "light" subexponential step-size distributions with

$$F(x+\Delta) = p(x) e^{-cx\log^{-\beta} x},$$

where $\beta > 0$, $c > 0$ and $p$ is $O$-regularly varying. On setting $R(x) = cx\log^{-\beta} x$ and noting that $yR'(y) = R(y) - \beta R(y)/\log y$, we find with Lemma B.2 that $x \mapsto F(x+\Delta)$ belongs to $\mathcal{S}d$. The small-step sequence found in the next proposition is not the best possible, but it suffices for our purposes.

PROPOSITION 8.4. *Suppose that $\mathsf{E}\{\xi\} = 0$ and $\mathsf{E}\{\xi^2\} = 1$, and consider the above setup.*

*The sequence $\{h_n \equiv \sqrt{n}\}$ is a truncation sequence, and $\{J_n \equiv \exp((c+\epsilon)^{1/\beta} n^{1/(2\beta)})\}$ is a corresponding small-steps sequence for any $\epsilon > 0$.*

PROOF. We only consider the global case, since the local case is similar. The arguments in the proof of Proposition 8.2 yield that $\{h_n\}$ is a truncation sequence. To show that $\{J_n\}$ is a corresponding small-steps sequence, we note that with Lemma 7.1(i), for $x \geq \exp((c+\epsilon)^{1/\beta} n^{1/(2\beta)})$,

$$\begin{aligned} &\frac{\mathsf{P}\{S_n > x, \xi_1 \leq \sqrt{n}, \ldots, \xi_n \leq \sqrt{n}\}}{n\overline{F}(x)} \\ &\quad \leq O(n^{-1}) \exp(-n^{-1/2} x - \log F(x+\Delta)) \\ &\quad \leq O(n^{-1}) \exp(-x[n^{-1/2} - (c+\epsilon/2)\log^{-\beta} x]), \end{aligned}$$

which is $o(1)$ since $\log^{-\beta}(x) \leq (c+\epsilon)^{-1} n^{-1/2}$. $\square$

When $p$ is intermediate regularly varying, we find an insensitivity sequence as in the previous two subsections, so that the next theorem follows from Theorem 2.1. To the best of our knowledge, the theorem is the first large-deviation result for (special cases of) the family under consideration.

THEOREM 8.4. *Let the assumptions of Proposition 8.4 hold, and suppose that $p$ is intermediate regularly varying at infinity.*

*For any $\{x_n\}$ with $x_n \gg n^{1/(2\beta)}$, we have $\mathsf{P}\{S_n \in x+\Delta\} \sim nF(x+\Delta)$ uniformly for $x \geq \exp(x_n)$.*

**9. Examples with infinite variance.** It is the aim of this section to work out our main theorem for classes of step-size distributions with infinite variance. Karamata's theory of regular variation and its ramifications provide the required additional structure.

9.1. *Infinite variance and a heavy right tail.* Having investigated the case where $F$ is attracted to a normal distribution, it is natural to also study the complementary case. We work within the framework of Karamata theory; see Appendix A.

We need three assumptions. Our main assumption is that

$$\overline{G}(x) \asymp x^{-2}\mu_2(x). \tag{25}$$

It is a well-known result due to Lévy that the "lower bound" part ensures that $F$ does not belong to the domain of partial attraction of the normal distribution. For more details we refer to Maller [28, 29]. Note that the "upper bound" part is exactly (21); it is shown by Feller [16] that this is equivalent with the existence of sequences $\{E_n\}$ and $\{F_n\}$ such that every subsequence of $\{S_n/E_n - F_n\}$ contains a further subsequence, say $\{n_k\}$, for which $S_{n_k}/E_{n_k} - F_{n_k}$ converges in distribution to a nondegenerate random variable. In that case, $\{S_n/E_n - F_n\}$ is called *stochastically compact*. Note that the required nondegeneracy is the only difference with $\{S_n/E_n - F_n\}$ being stochastically bounded; further details can, for instance, be found in Jain and Orey [24].

When interpreting (25), it is important to realize the following well-known fact. If $F$ is attracted to a stable law with index $\alpha \in (0,2)$, then the tails must be regularly varying, and application of Karamata's theorem shows that $\alpha\overline{G}(x) \sim (2-\alpha)x^{-2}\mu_2(x)$. Therefore, our assumption (25) is significantly more general.

Our second assumption is that the left tail of $F$ is not heavier than the right tail:

$$\limsup_{x\to\infty} \frac{\overline{G}(x)}{\overline{F}(x)} < \infty. \tag{26}$$

In the next subsection, we investigate the complementary case with a heavier left tail.

Our third assumption, which is formulated in terms of the $a_n$ defined in Section 7, ensures that $F$ is sufficiently centered:

$$\limsup_{n\to\infty} \frac{n|\mu_1(a_n)|}{a_n} < \infty. \tag{27}$$

This assumption often follows from (25), as shown in the next lemma. The lemma also records other important consequences of (25), and relies completely on the seminal work on $O$-regular variation by Bingham, Goldie and Teugels [3]. Item (i) is due to Feller [16], but the reader is advised to also refer to the extended and corrected treatment in [3].

LEMMA 9.1. *Equation (25) is equivalent to the following two statements:*

(i) *$\mu_2$ is $O$-regularly varying with Matuszewska indices $0 \le \beta_{\mu_2} \le \alpha_{\mu_2} < 2$.*

(ii) *$\overline{G}$ is $O$-regularly varying with Matuszewska indices $-2 < \beta_{\overline{G}} \le \alpha_{\overline{G}} \le 0$.*

*Moreover, under (25), we automatically have (27) if either $\beta_{\overline{G}} > -1$, or if $\alpha_{\overline{G}} < -1$ and $\mathsf{E}\{\xi\} = 0$.*

PROOF. All cited theorems in this proof refer to [3]. The equivalence of (25) and (i), (ii) follows from Theorem 2.6.8(c) and Theorem 2.6.8(d). If $\beta_{\overline{G}} > -1$, then we have $\limsup_{x\to\infty} x^{-1}\int_0^x yG(dy)/\overline{G}(x) < \infty$ by Theorem 2.6.8(d). Similarly, if $\mathsf{E}\{|\xi|\} < \infty$ and $\alpha_{\overline{G}} < -1$, then $\limsup_{x\to\infty} x^{-1}\int_x^\infty yG(dy)/\overline{G}(x) < \infty$ by Theorem 2.6.7(a), (c). $\square$

The next proposition gives appropriate truncation and small-steps sequences.

PROPOSITION 9.1. *Suppose that (25), (26) and (27) hold. Moreover, if $T < \infty$, suppose that $x \mapsto F(x+\Delta)$ is $O$-regularly varying with upper Matuszewska index strictly smaller than $-1$.*

*Given some $\{t_n\}$ with $n\overline{G}(t_n) = o(1)$, there exists some $\gamma > 0$ such that, with*

$$h_n \equiv \frac{t_n}{-2\gamma \log[n\overline{G}(t_n/2)]},$$

*$\{h_n\}$ is a truncation sequence. Moreover, $\{J_n \equiv t_n/2\}$ is then an $h$-small-steps sequence.*

PROOF. We first show that $\{h_n\}$ is a truncation sequence. Our assumption on $F(x+\Delta)$ guarantees that it is almost decreasing. In view of Lemmas 6.1 and 9.1, it suffices to show that $n\overline{F}(h_n) = o(1)$. The first step is to prove that $h_n \to \infty$, for which we use the bound $\overline{G}(x) \geq x^{-2}$ for large $x$ (see Theorem 2.2.7 in [3]): we have that

$$h_n \geq \frac{t_n}{-2\gamma \log[4nt_n^{-2}]} \geq \frac{t_n}{-2\gamma \log(n) + 4\gamma \log(t_n/2)} \geq \frac{t_n}{4\gamma \log(t_n/2)},$$

which exceeds any given number for large $n$. Relying on the fact that $h_n \to \infty$, the Potter-type bounds of Proposition 2.2.1 in [3] yield that for some $C' > 0$, provided $n$ is large, $\overline{G}(t_n/2)/\overline{G}(h_n) \geq C'(t_n/(2h_n))^{-2}$. Hence, by definition of $h_n$, as $n \to \infty$,

$$n\overline{G}(h_n) \leq (C')^{-1}(-2\gamma \log[n\overline{G}(t_n/2)])^2 n\overline{G}(t_n/2) = o(1).$$

This proves in particular that $n\overline{F}(h_n) = o(1)$, so that $\{h_n\}$ is a truncation sequence.

We now prove that $\{t_n/2\}$ is a small-steps sequence. The idea is to apply Lemma 7.1(ii), for which we need $h_n \geq a_n$. In fact, we have $h_n \gg a_n$; this follows from the fact that $n\overline{G}(a_n)$ is bounded away from zero [note that $\overline{G}(x) \asymp Q(x)$ by (25)] in conjunction with our observation that $n\overline{G}(h_n) = o(1)$. Throughout the proof, let $C < \infty$ be a generic constant which can change from line to line.

First consider the global case $T = \infty$. Lemma 7.1(ii) shows that for any $x \geq 0$,

$$\sup_{z \geq x} \mathsf{P}\{S_n > z, \xi_1 \leq h_n, \dots, \xi_n \leq h_n\} \leq C \exp(-x/h_n).$$

This shows that for $\gamma > 2$, by (26), the aforementioned Potter-type bound and the definition of $h_n$,

$$\begin{aligned}
&\sup_{x \geq J_n} \sup_{z \geq x} \frac{\mathsf{P}\{S_n > z, \xi_1 \leq h_n, \dots, \xi_n \leq h_n\}}{n\overline{F}(x)} \\
&\quad \leq C \sup_{x \geq 1} \frac{e^{-(t_n/(2h_n))x}}{n\overline{G}(xt_n/2)} \leq C \sup_{x \geq 1} x^2 e^{(1/2)\gamma \log[n\overline{G}(t_n/2)]x} \frac{\exp(-(t_n/(4h_n)x)}{n\overline{G}(t_n/2)} \\
&\quad \leq C \sup_{x \geq 1} x^2 e^{-x} \frac{\exp(-(t_n/(4h_n))}{n\overline{G}(t_n/2)} \leq C(n\overline{G}(t_n/2))^{\gamma/2-1} = o(1).
\end{aligned}$$

Similar ideas are used to prove the local case, but now we also need the concentration result of Proposition 7.1. Since $h_n \gg a_n$, we use this proposition in conjunction with Lemma 7.1(ii) to conclude that

$$\sup_{z \geq x} \mathsf{P}\{S_n \in z + \Delta, \xi_1 \leq h_n, \dots, \xi_n \leq h_n\} \leq C a_n^{-1} \exp(-x/h_n).$$

To prove the proposition, by (26) it therefore suffices to show that for some $\gamma > 0$,

$$(n\overline{F}(t_n/2))^\gamma = o(na_n F(t_n/2 + \Delta)).$$

The assumption on $F(x+\Delta)$ is equivalent with $\overline{F}(x) \asymp xF(x+\Delta)$ by Corollary 2.6.4 of [3]. Therefore, it is enough to prove the above equality with $F(t_n/2+\Delta)$ replaced by $t_n^{-1}\overline{F}(t_n/2)$.

On combining the assumption on $F(x+\Delta)$ with (26), we obtain $\overline{G}(x) \asymp \overline{F}(x) \asymp xF(x+\Delta)$. Hence, $\overline{G}$ has bounded decrease, which implies (see Proposition 2.2.1 of [3]) that there exists some $\eta > 0$ such that

$$\begin{aligned}\frac{t_n}{a_n}[n\overline{F}(t_n/2)]^{\gamma-1} &\le \frac{t_n}{a_n}[n\overline{G}(t_n/2)]^{\gamma-1} \le C\frac{t_n}{a_n}\left(\left[\frac{t_n}{a_n}\right]^{-\eta} n\overline{G}(a_n)\right)^{\gamma-1} \\ &\le C\frac{t_n}{a_n}\left(\left[\frac{t_n}{a_n}\right]^{-\eta} nQ(a_n)\right)^{\gamma-1}.\end{aligned}$$

This upper bound vanishes if we choose $\gamma > 1 + 1/\eta$. □

Let us now define $b_n \equiv h_n$. Since $\{S_n/a_n\}$ is tight under the assumptions of the preceding proposition (see, e.g., [24], Proposition 1.2), and since we have shown in its proof that $b_n \gg a_n$, we immediately conclude that $S_n/b_n$ converges in distribution to zero. In particular, $\{b_n\}$ is a natural-scale sequence.

It remains to choose a corresponding insensitivity sequence. This can immediately be done under the assumption that $x \mapsto F(x+\Delta)$ is intermediate regularly varying (see Appendix A). Indeed, since $b_n \ll t_n/2$, we may set $I_n = t_n/2$ and conclude with Corollary 2.2I of [8] that $\{I_n\}$ is an insensitivity sequence.

We have proved that the next theorem follows from Theorem 2.1. The theorem has a long history. In the global regularly varying case, it is due to Heyde [21]; S. Nagaev [34] ascribes it to Tkachuk. For a recent account, see Borovkov and Boxma [6]. Heyde [20] studies the nonregularly varying case, but only proves the right order of $\mathsf{P}\{S_n > x\}$; related results have been obtained by Cline and Hsing [9]. In the local case, only the regularly varying case has been investigated. Our theorem then reproduces the large-deviation theorem in Doney [12] in the infinite-mean case, while significantly improving upon the results in Doney [11] in the complementary case.

THEOREM 9.1. *Let the assumptions of Proposition 9.1 hold, and suppose that $x \mapsto F(x+\Delta)$ is intermediate regularly varying at infinity.*

*For any $\{x_n\}$ with $n\overline{F}(x_n) = o(1)$, we have $\mathsf{P}\{S_n \in x+\Delta\} \sim nF(x+\Delta)$ uniformly for $x \ge x_n$.*

9.2. *Finite mean, infinite variance, and a heavy left tail.* In this subsection, we investigate the case when the left tail is heavier than the right tail, and this tail causes $\xi$ to be integrable yet to have an infinite second moment. It is our aim to recover the big-jump result derived by Rozovskii [41] in this context, and to extend it to the local case.

More precisely, we assume that:

- $x \mapsto F(-x)$ is regularly varying at infinity with index $-\alpha$ for some $\alpha \in (1,2)$,
- $x \mapsto \overline{F}(x)$ is regularly varying at infinity with index $-\beta$ for some $\beta > \alpha$, and
- $\mathsf{E}\{\xi\} = 0$.

Under these assumptions, $F$ belongs to the domain of attraction of the $\alpha$-stable law with a Lévy measure that vanishes on the positive half-line. The theory on domains of attraction (e.g., [17], Section XVII.5) immediately implies that $\{b_n\}$ determined by $\Gamma(3-\alpha)n\mu_2(b_n) = (\alpha-1)b_n^2$ is a natural-scale sequence. Note that this sequence is regularly varying with index $1/\alpha$, and that $n\overline{G}(b_n)$ tends to a constant. The next proposition shows how $\{h_n\}$ and $\{J_n\}$ can be chosen under a condition which should be compared with [41], (1.19).

PROPOSITION 9.2. *Suppose that the above three assumptions hold, and that*

$$\limsup_{n\to\infty} \frac{F(-b_n/[\log n]^{1/\alpha})}{(\log n)F(-b_n)} \le 1. \tag{28}$$

*Furthermore, if $T < \infty$, suppose that $F(x+\Delta)$ is regularly varying.*

*The sequence $\{h_n \equiv (\frac{\beta-\alpha}{\alpha-1}\log n)^{-1/\alpha} b_n\}$ is a truncation sequence. Moreover, given $t > 1$, if we set*

$$J_n = t\left(\frac{\beta-\alpha}{\alpha-1}\log n\right)^{(\alpha-1)/\alpha} b_n,$$

*then $\{J_n\}$ is an $h$-small-steps sequence.*

PROOF. To see that $\{h_n\}$ is a truncation sequence, we use Lemma 6.1 and the elementary bounds

$$\begin{aligned} n\overline{F}(h_n) &\le n h_n^{-3\beta/4-\alpha/4} \le h_n^{-(\beta-\alpha)/2} h_n^{-\beta/4-3\alpha/4} \le h_n^{-(\beta-\alpha)/2} nF(-h_n) \\ &\le 2(\log n) h_n^{-(\beta-\alpha)/2} nF(-b_n) \le 4(\log n) h_n^{-(\beta-\alpha)/2}, \end{aligned}$$

where we have used (28). Since $\{h_n\}$ is regularly varying with index $1/\alpha$, this upper bound tends to zero.

We next concentrate on $\{J_n\}$, for which we use Lemma 7.1(iii). Choose $0 < 4\epsilon < t-1$. If $\int_0^\infty u^2 F(du) = \infty$, application of Karamata's theorem (on the right tail) shows that

$$h_n^{-2} \int_0^{h_n} u^2 F(du) = (1+o(1))\overline{F}(h_n) = o(F(-h_n)).$$

We reach the same conclusion in the complementary case $\int_0^\infty u^2 F(du) < \infty$. Using (28) we obtain that, for large $n$,

$$\begin{aligned}
&\frac{n}{h_n^2}\int_0^{h_n} u^2 F(du) + (1+\epsilon)\frac{\Gamma(2-\alpha)}{\alpha-1} nF(-h_n) \\
&\quad \le (1+2\epsilon)\frac{\Gamma(2-\alpha)}{\alpha-1} nF(-h_n) \le (1+3\epsilon)\frac{1}{\alpha}\frac{F(-h_n)}{F(-b_n)} \le t\frac{\beta-\alpha}{\alpha(\alpha-1)}\log n. \qquad (29)
\end{aligned}$$

We now have all the prerequisites to prove the claim in the global case, that is, for $T = \infty$. Indeed, we need to show that, for the $\{h_n\}$ and $\{J_n\}$ given above,

$$\sup_{x \ge J_n}\left[-\frac{x}{h_n} + t\frac{\beta-\alpha}{\alpha(\alpha-1)}\log n - \log n - \log \overline{F}(x)\right] \to -\infty.$$

Fix some $0 < \eta < (t-1)(\beta-\alpha)$. The elementary estimate $\overline{F}(x) \ge x^{-\beta-\eta}$ (for large $x$) yields an upper bound for which the supremum is attained at $J_n$ for large $n$. We conclude that the left-hand side of the preceding display is bounded from above by

$$\begin{aligned}
&-\frac{J_n}{h_n} + t\frac{\beta-\alpha}{\alpha(\alpha-1)}\log n - \log n + \left(\frac{\beta+\eta}{\alpha}\right)\log n \\
&\quad = -(t-1)\frac{\beta-\alpha}{\alpha}\log n + \frac{\eta}{\alpha}\log n \to -\infty.
\end{aligned}$$

It remains to treat the local case $T < \infty$, for which we use similar arguments based on Chebyshev's inequality. The bound (29), in conjunction with Proposition 7.1(ii) and the fact that $h_n \le b_n$, shows that it suffices to prove

$$\sup_{x \ge J_n}\left[-\frac{x}{h_n} + t\frac{\beta-\alpha}{\alpha(\alpha-1)}\log n - \log n - \log F(x+\Delta) - \log h_n\right] \to -\infty.$$

The index of regular variation of $x \mapsto F(x+\Delta)$ is necessarily $-\beta-1$ by Karamata's theorem. We can now repeat the reasoning for the global case, observing that $-\log h_n + \log J_n = o(\log J_n)$. $\square$

To gain some intuition for the above proposition, it is instructive to see how $\{h_n\}$ and $\{J_n\}$ arise as a result of an optimization procedure similar to the finite-variance heuristic given at the end of Section 7. Suppose for simplicity that $F(-x) = x^{-\alpha}$ and that $1 + o(1)$ may be read as 1. The last but one bound in (29) shows that $J_n$ must exceed $b_n^\alpha h_n^{-\alpha+1} - h_n \log n + \beta/\alpha h_n \log n$. Now optimize this bound with respect to $h_n$ to find the sequences of the proposition.

We also remark that our reasoning immediately allows for a relaxation of the assumptions on the right tail, for instance in terms of $O$-regular variation. In fact, the proof shows that Karamata assumptions on the right tail can be avoided altogether by assuming that $\int_0^\infty u^2 F(du) < \infty$, and then replacing $\beta$ in the statement by $\inf\{\gamma : \liminf_{x\to\infty} x^\gamma \overline{F}(x) > 0\}$. Still, regular variation of the left tail is essential in order to apply Lemma 7.1(iii), which relies on a Tauberian argument.

The next theorem is a corollary of the preceding proposition in conjunction with Theorem 2.1. In the global case it has been obtained by Rozovskii [41], Corollary 2A.

Theorem 9.2. *Let the assumptions of Proposition 9.2 hold. For any $t > 1$, we have $\mathsf{P}\{S_n \in x + \Delta\} \sim nF(x + \Delta)$ uniformly for $x \geq t(\frac{\beta-\alpha}{\alpha-1} \log n)^{(\alpha-1)/\alpha} b_n$.*

9.3. *Infinite mean and a heavy left tail.* In this subsection we consider the case when the left tail is heavier than the right tail, and when $\xi$ fails to be integrable. This situation has recently been studied by Borovkov [5]; we include it here to show an interesting contrast with the preceding subsection, which is perhaps surprising in view of the unified treatment in Section 9.1 for balanced tails.

We assume that:

- $x \mapsto F(-x)$ is regularly varying at infinity with index $-\alpha$ for some $\alpha \in (0, 1)$, and
- $x \mapsto \overline{F}(x)$ is regularly varying at infinity with index $-\beta$ for some $\beta > \alpha$.

Under these assumptions, $F$ is in the domain of attraction of the unbalanced $\alpha$-stable law, and $\{b_n\}$ with $b_n = \inf\{x : F(-x) < 1/n\}$ is a natural-scale sequence.

The following proposition shows that, in the present situation, one can take a small-steps sequence which is fundamentally different from the one in Section 9.2.

Proposition 9.3. *Suppose that the above two assumptions hold. If $T < \infty$, also suppose that $F(x + \Delta)$ is regularly varying.*

*The sequence $\{h_n \equiv n^{1/\beta}\}$ is a truncation sequence. Moreover, for any given $\epsilon > 0$, the sequence $\{J_n \equiv n^{1/\beta+\epsilon}\}$ is an $h$-small-steps sequence.*

PROOF. The proof is modeled after the proof of Proposition 9.2. It becomes clear with Lemma 6.1 that $\{h_n\}$ is a natural-scale sequence.

We next apply Lemma 7.1(iv). If $\int_0^\infty uF(du) = \infty$, we apply Karamata's theorem and see that $nh_n^{-1}\int_0^{h_n} uF(du)$ is $o(nF(-h_n))$; otherwise we conclude this immediately. Similarly, $nh_n^{-2}\int_0^{h_n} u^2F(du)$ is always $o(nF(-h_n))$. This shows that, for sufficiently large $n$, $n\log\int_{-\infty}^{h_n} e^{u/h_n}F(du) \le 0$. Therefore, if $T = \infty$, it suffices to observe that $h_n, J_n$ satisfy

$$\lim_{n\to\infty}\sup_{x\ge J_n}\frac{\exp(-x/h_n)}{nx^{-\beta-\epsilon}} = 0.$$

The local case is similar. □

The next theorem, which is new in the local case, immediately follows from the preceding proposition in conjunction with Theorem 2.1.

THEOREM 9.3. *Let the assumptions of Proposition 9.3 hold. For any* $\{x_n\}$ *with* $nF(-x_n) = o(1)$, *we have* $\mathsf{P}\{S_n \in x+\Delta\} \sim nF(x+\Delta)$ *uniformly for* $x \ge x_n$.

## APPENDIX A: SOME NOTIONS FROM KARAMATA THEORY

We recall some useful notions from Karamata theory for the reader's convenience. A positive, measurable function $f$ defined on some neighborhood of infinity is *$O$-regularly varying* (at infinity) if

$$0 < \liminf_{x\to\infty}\frac{f(xy)}{f(x)} \le \limsup_{x\to\infty}\frac{f(xy)}{f(x)} < \infty.$$

This is equivalent to the existence of some (finite) $\alpha_f, \beta_f$ with the following properties. For any $\alpha > \alpha_f$, there exists some $C = C(\alpha)$ such that for any $Y > 1$, $f(xy)/f(x) \le C(1+o(1))y^\alpha$ uniformly in $y \in [1, Y]$. Similarly, for any $\beta < \beta_f$, there exists some $D = D(\beta)$ such that for any $Y > 1$, $f(xy)/f(x) \ge D(1+o(1))y^\beta$ uniformly in $y \in [1, Y]$. The numbers $\alpha_f$ and $\beta_f$ are called the *upper* and *lower Matuszewska indices*, respectively. We refer to [3], Chapter 2 for more details.

A positive, measurable function $f$ defined on some neighborhood of infinity is *intermediate regularly varying* (at infinity) if

$$\lim_{y\downarrow 1}\liminf_{x\to\infty}\frac{f(xy)}{f(x)} = \lim_{y\downarrow 1}\limsup_{x\to\infty}\frac{f(xy)}{f(x)} = 1.$$

Intermediate regular variation has been introduced by Cline [8]. Cline also shows that an intermediate regularly varying function is necessarily $O$-regularly varying. Note that regular variation implies intermediate regular variation.

## APPENDIX B: THE CLASS $\mathcal{S}D$ OF SUBEXPONENTIAL DENSITIES

We say that a function $H:\mathbb{R}\to\mathbb{R}_+$ belongs to the class $\mathcal{S}d$ if $H\in\mathcal{L}$ and

$$\lim_{x\to\infty}\frac{\int_0^{x/2}H(y)H(x-y)\,dy}{H(x)}=\int_0^\infty H(y)\,dy<\infty.$$

It is important to realize that it is possible to determine whether $H$ belongs to $\mathcal{S}d$ by considering its restriction to the positive half-line. Under the extra assumptions that $H$ be monotone and supported on the positive half-line, the requirement $H\in\mathcal{L}$ is redundant and the class is usually referred to as $\mathcal{S}^*$.

This section aims to present criteria for assessing whether a function $H\in\mathcal{L}$ of the form

$$H(x)=p(x)e^{-R(x)} \tag{30}$$

belongs to $\mathcal{S}d$, where $p$ is $O$-regularly varying.

LEMMA B.1. *Consider $H\in\mathcal{L}$ of the form (30), where $p$ is $O$-regularly varying. Suppose that there exists an eventually concave function $z\ge 0$ such that $\limsup xz'(x)/z(x)<1$ and the function $R(x)/z(x)$ is eventually nonincreasing. If, moreover, $R(x)\gg\log x$, then we have $H\in\mathcal{S}d$.*

PROOF. It follows from $H\in\mathcal{L}$ that there is some $h$ with $h(x)\le x/2$, $h(x)\to\infty$, and $H(x-y)\sim H(x)$ uniformly for $y\le h(x)$. Therefore, we have

$$\begin{aligned}\int_0^{h(x)}H(y)H(x-y)\,dy&\sim H(x)\int_0^{h(x)}H(y)\,dy\\&\sim H(x)\int_0^\infty H(y)\,dy.\end{aligned}$$

It therefore suffices to show that the integral over the interval $(h(x),x/2]$ is $o(H(x))$. Exploiting the assumptions on $R$ and $z$, the proof of Theorem 2 of Shneer [43] in conjunction with Property 2 in [42] shows that there exists an $\alpha\in(0,1)$ such that $R(x)-R(x-y)\le\alpha yR(x)/x$ for $0\le y\le x/2$ for large $x$. Moreover, since $x\mapsto R(x)/x$ is ultimately nonincreasing, we have $R(x)-R(x-y)-R(y)\le(\alpha-1)R(y)$ for $h(x)\le y\le x/2$. The imposed $O$-regular variation of $p$ implies $\sup_{u\in[1/2,1]}p(ux)/p(x)=O(1)$ and $p(x)\le x^\eta$ for some $\eta<\infty$ and large enough $x$, showing that

$$\begin{aligned}\frac{\int_{h(x)}^{x/2}H(y)H(x-y)\,dy}{H(x)}&\le\int_{h(x)}^{x/2}\frac{p(y)p(x-y)}{p(x)}e^{-(1-\alpha)R(y)}\,dy\\&\le O(1)\int_{h(x)}^{x/2}p(y)e^{-(1-\alpha)R(y)}\,dy\\&\le O(1)\int_{h(x)}^{x/2}y^{-2}\,dy,\end{aligned}$$

where we have also used $R(x) \gg \log x$ to obtain the last inequality. □

The next lemma is inspired by Theorem 3.6(b) of Klüppelberg [26].

LEMMA B.2. *Consider $H \in \mathcal{L}$ of the form (30), where $p$ is $O$-regularly varying. Suppose that $R$ is differentiable and that $R'$ is ultimately nonincreasing. If $\int_M^\infty e^{yR'(y)} H(y)\,dy < \infty$ for some $M < \infty$, then $H \in \mathcal{S}d$.*

PROOF. As in the proof of the previous lemma, it suffices to bound $H(y) \times H(x-y)/H(x)$ for $y \in (h(x), x/2]$. We have $x - y \geq y$ for $y \leq x/2$, implying that

$$R(x) - R(x-y) \leq yR'(x-y) \leq yR'(y).$$

Note that $p(x-y)/p(x) = O(1)$ uniformly for $y \leq x/2$ since $p$ is $O$-regular varying, yielding

$$\frac{\int_{h(x)}^{x/2} H(y)H(x-y)\,dy}{H(x)} \leq O(1) \int_{h(x)}^{\infty} e^{yR'(y)} H(y)\,dy,$$

which vanishes by assumption. □

**Acknowledgments.** The authors are grateful to Thomas Mikosch for several stimulating and helpful discussions, and to Ken Duffy for providing detailed comments. A. B. Dieker's research has been carried out in part when he was with University College Cork, Cork, Ireland, and CWI, Amsterdam, The Netherlands. The research was conducted in part when V. Shneer was with Heriot-Watt University and when D. Denisov was with EURANDOM.

D. Denisov
Heriot-Watt University
School of Mathematical and Computer Sciences
Edinburgh EH14 4AS
United Kingdom
E-mail: denisov@ma.hw.ac.uk

A. B. Dieker
IBM Research
T. J. Watson Research Center
P.O. Box 218
Yorktown Heights, New York 10598
USA
E-mail: tondieker@us.ibm.com

V. Shneer
EURANDOM
P.O. Box 513
5600 MB Eindhoven
The Netherlands
E-mail: shneer@eurandom.tue.nl